\documentclass[12pt, reqno]{amsart}
\usepackage{amsmath, amstext, amsbsy, amssymb, amscd}

\newtheorem{theorem}{Theorem}[section]
\newtheorem{lemma}[theorem]{Lemma}
\newtheorem{proposition}[theorem]{Proposition}

\theoremstyle{definition}
\newtheorem{definition}[theorem]{Definition}

\theoremstyle{remark}

\numberwithin{equation}{section}

\newcommand{\C}{ \mathbb C }

\newcommand{\fM}{\frak M}
\newcommand{\cE}{\mathcal E}

\newcommand{\Ext}{\text{\rm Ext}}

\newcommand{\Grass}{\text {\rm Grass}}

\newcommand{\Hilb}{\text{\rm Hilb}}

\newcommand{\M}{\mathfrak M}

\newcommand{\pd}{\text{\rm PD}}
\newcommand{\PP}{{\mathbb P}}
\newcommand{\rev}{\widetilde {ev}_1}
\newcommand{\rf}{\tilde {f}_{1, 0}}
\newcommand{\SL}{\text{\rm SL}}
\newcommand{\Supp}{\text{\rm Supp}}

\newcommand{\vac}{|0\rangle}

\newcommand{\Xn}{ X^{[n]}}

\newcommand{\Z}{ \mathbb Z }

{\vskip-\lastskip\medskip
  \noindent
  {\em #1.}\enspace
  }%
{\qed\par\medskip
  }

\begin{document}
\title[Gromov-Witten invariants of Hilbert schemes]
{Gromov-Witten invariants \\
of the Hilbert scheme of $3$-points on ${\mathbb P}^2$}

\author[Dan Edidin]{Dan Edidin$^1$}
\address{Department of Mathematics, University of Missouri, Columbia, MO
65211, USA} \email{edidin@math.missouri.edu}
\thanks{${}^1$Partially supported by an NSF grant and an NSA grant}
\author[Wei-Ping Li]{Wei-Ping Li$^2$}
\address{Department of Mathematics, HKUST, Clear Water Bay, Kowloon, Hong
Kong } \email{mawpli@ust.hk}
\thanks{${}^2$Partially supported by the grant HKUST6170/99P}
\author[Zhenbo Qin]{Zhenbo Qin$^3$}
\address{Department of Mathematics, University of Missouri, Columbia, MO
65211, USA} \email{zq@math.missouri.edu}
\thanks{${}^3$Partially supported by an NSF grant and
a University of Missouri Research Board grant}


\begin{abstract}
Using obstruction bundles, composition law and localization formula,
we compute certain $3$-point genus-$0$ Gromov-Witten invariants of
the Hilbert scheme of $3$-points on the complex projective plane.
Our results partially verify Ruan's conjecture
about quantum corrections for this Hilbert scheme.
\end{abstract}

\maketitle
\date{}
%\tableofcontents

%%
%%
%%
%%
%%
%%
%%
\section{\bf Introduction}

Motivated by the pioneering work of Nakajima and Grojnowski \cite{Nak, Gro},
there have been intensive studies of the cohomology ring structure of
the Hilbert schemes of points on a smooth algebraic surface
(e.g. \cite{Leh, L-S, LQW1, LQW2, LQW3, Q-W, Go2}). While our understanding of
this ordinary cohomology ring structure has deepened rapidly,
the quantum cohomology ring structure of
these Hilbert schemes remains to be a mystery. A limited progress
to the quantum cohomology ring structure has been made in \cite{L-Q}
where certain $1$-point genus-$0$ Gromov-Witten invariants of
these Hilbert schemes have been determined. These $1$-point invariants
come from the contributions of curves contracted by the Hilbert-Chow map
from the Hilbert schemes to the symmetric products of the surface.

In this paper, we study $3$-point genus-$0$ Gromov-Witten invariants of
the Hilbert scheme $({\mathbb P}^2)^{[3]}$ of $3$-points on
the complex projective plane ${\mathbb P}^2$. Again, we are primarily
interested in those invariants which come from the contributions of
curves contracted by the Hilbert-Chow map (\ref{H_C}).
These curves are homologous to $d \beta_3$ for some positive integer $d$,
where $\beta_3 \subset ({\mathbb P}^2)^{[3]}$ is
the rational curve defined by
\begin{eqnarray*}
\beta_3 = \{ \xi + x_2| \, \ell(\xi) = 2, \,\, \Supp(\xi) = x_1 \}
\end{eqnarray*}
with $x_1$ and $x_2$ being two fixed distinct points of the projective plane
$X = \PP^2$.

To state our main results, we introduce some notations.
Let $H^*(X^{[3]})$ and $H_*(X^{[3]})$ be the cohomology and homology of
$X^{[3]}$ with $\C$-coefficients. For $i = 2, 4, 6, 8, 10$,
a linear basis $\frak B_i$ of $H_i(X^{[3]})$ in terms of
the Heisenberg operators introduced in \cite{Nak, Gro}
can be determined (see Lemma~\ref{lemmacohoclass}
and Definition~\ref{def_B_i} for details).
For $\alpha_1, \ldots, \alpha_k \in H^*(X^{[3]})$,
we use $\langle \alpha_1, \ldots, \alpha_k \rangle_{0, d}$ to
stand for the $k$-point genus-$0$ Gromov-Witten invariant
$\langle \alpha_1, \ldots, \alpha_k \rangle_{0, d \beta_3}$.
Now the $3$-point genus-$0$ Gromov-Witten invariants
$\langle \alpha_1, \alpha_2, \alpha_3 \rangle_{0, d}$ of $X^{[3]}$
are reduced either to the $2$-point invariants
$\langle \pd(A_1), \pd(A_2) \rangle_{0, d}$
with $A_1 \in \frak B_6$ and $A_2 \in \frak B_8$,
or to the $3$-point invariants
$\langle \pd(A_1), \pd(A_2), \pd(A_3) \rangle_{0, d}$
with $A_1, A_2, A_3 \in \frak B_8$.
Here PD denotes the Poincar\' e duality.
Our main results are the following.

\begin{theorem}  \label{intro_thm_2_point}
Let $X = \PP^2$, and $\frak B_6$ and $\frak B_8$ be
defined in Definition~\ref{def_B_i}.
Let $d \ge 1$, $A_1 \in \frak B_6$ and $A_2 \in \frak B_8$.
Let $x, \ell$ be a point and a line in $X$ respectively.
Then, $\langle \pd(A_1), \pd(A_2) \rangle_{0, d}$ is zero unless
the pair $(A_1, A_2)$ is one of the following:
\par
{\rm (i)}
$(\mathfrak a_{-2}(X)\mathfrak a_{-1}(x) \vac,
\mathfrak a_{-1}(X)\mathfrak a_{-2}(\ell)\vac)$
\par
{\rm (ii)} $(\mathfrak a_{-2}(\ell)\mathfrak a_{-1}(\ell)\vac,
\mathfrak a_{-2}(X)\mathfrak a_{-1}(\ell)\vac)$
\par
{\rm (iii)} $(\mathfrak a_{-3}(\ell)\vac, \mathfrak a_{-3}(X)\vac)$.
\par\noindent
Moreover, $\langle \pd(A_1), \pd(A_2) \rangle_{0, d} = 12/d$
in cases (i) and (ii).
\end{theorem}

\begin{theorem}  \label{intro_thm_3_point}
Let $X = \PP^2$, and $\frak B_8$ be defined in Definition~\ref{def_B_i}.
Let $\ell \subset X$ be a line.
Let $d \ge 1$, $f(d) = d \, \langle \pd(\frak a_{-3}(\ell)\vac),
\pd(\frak a_{-3}(X)\vac) \rangle_{0, d}$,
and $A_1, A_2, A_3 \in \frak B_8$.
Then, the $3$-point genus-$0$ Gromov-Witten invariant
$\langle \pd(A_1), \pd(A_2), \pd(A_3) \rangle_{0, d}$ is zero unless
the unordered triple $(A_1, A_2, A_3)$ is one of the following:
\par
{\rm (i)}
$(\frak a_{-2}(X)\frak a_{-1}(\ell)\vac,
  \frak a_{-2}(X)\frak a_{-1}(\ell)\vac,
  \frak a_{-1}(X)\frak a_{-2}(\ell)\vac)$
\par
{\rm (ii)} $(\mathfrak a_{-3}(X)\vac, \mathfrak a_{-3}(X)\vac,
\frak a_{-2}(X)\frak a_{-1}(\ell)\vac)$
\par
{\rm (iii)} $(\mathfrak a_{-3}(X)\vac, \mathfrak a_{-3}(X)\vac,
\frak a_{-1}(X)\frak a_{-2}(\ell)\vac)$
\par
{\rm (iv)} $(\mathfrak a_{-3}(X)\vac, \mathfrak a_{-3}(X)\vac,
\mathfrak a_{-3}(X)\vac)$.
\par\noindent
Moreover, $\langle \pd(A_1), \pd(A_2), \pd(A_3) \rangle_{0, d}
= -24$ for case (i); for cases (ii) and (iii),
$\langle \pd(A_1), \pd(A_2), \pd(A_3) \rangle_{0, d} = -2f(d)$;
for case (iv),
\begin{eqnarray*}
& &\langle \pd(A_1), \pd(A_2), \pd(A_3) \rangle_{0, d} \\
&=&-162-15f(d) +6\sum\limits_{0 < d_1 < d} f(d_1)
+{1 \over 3} \sum\limits_{0 < d_1 < d} f(d_1)f(d-d_1).
\end{eqnarray*}
\end{theorem}

These two theorems are proved by using obstruction bundles and
composition laws in Sect.~\ref{sect_GW},
which generalizes the earlier methods in \cite{L-Q}.
In view of our theorems, to compute all the $3$-point invariants
$\langle \alpha_1, \alpha_2, \alpha_3 \rangle_{0, d}$ of $X^{[3]}$,
it remains to determine the $2$-point invariant
$\langle \pd(\frak a_{-3}(\ell)\vac),
\pd(\frak a_{-3}(X)\vac) \rangle_{0, d}$. In Sect.\ref{sect_loc},
using the standard $(\C^*)^2$-action on $X = \PP^2$
and the virtual localization formula from \cite{G-P},
we reduce the computation of $\langle \pd(\frak a_{-3}(\ell)\vac),
\pd(\frak a_{-3}(X)\vac) \rangle_{0, d}$ to
a summation over stable graphs. Even though we could not simplify
this summation for a general $d$, we are able to calculate the summation
for $d \le 4$ by employing Mathematica.
This enables us to prove the following.

\begin{proposition} \label{intro_prop}
Let $X = \PP^2$, and $\ell \subset X$ be a line.
Then, the $2$-point genus-$0$ Gromov-Witten invariant
$\langle \pd(\frak a_{-3}(\ell)\vac),
\pd(\frak a_{-3}(X)\vac) \rangle_{0, d}$ is equal to
$-27$, $27/2$, $18$ and $27/4$
when $d$ is equal to $1$, $2$, $3$ and $4$ respectively.
\end{proposition}

One of our motivations for this present work is to verify
Ruan's conjecture in \cite{Ru2} about the quantum corrections
for crepant resolutions of orbifolds. The symmetric products of
a smooth projective surface are global orbifolds.
The Hilbert-Chow map (\ref{H_C}) presents the Hilbert schemes of
points on a smooth projective surface as crepant resolutions of
the symmetric products of the surface.
For the Hilbert scheme $({\mathbb P}^2)^{[3]}$, our results enable
us to verify Ruan's conjecture for those quantum corrections
not involving $\langle \pd(\frak a_{-3}(\ell)\vac),
\pd(\frak a_{-3}(X)\vac) \rangle_{0, d}$.
Since the verification involves only straight-forward computations,
we omit the details.

Finally, we remark that our methods
can be extended in several directions. First of all,
they can be used to compute many $3$-point Gromov-Witten invariants
of the Hilbert scheme $({\mathbb P}^2)^{[n]}$ for a general $n$.
Secondly, our methods of proving Theorem~\ref{intro_thm_2_point}
and Theorem~\ref{intro_thm_3_point} can be easily modified to work
for an arbitrary simply connected projective surface $X$.
In addition, the ideas of proving Proposition~\ref{intro_prop}
can be applied to other toric surfaces.
We leave the details to the interested readers.

%\medskip\noindent
%{\bf Conventions:}
%Unless otherwise specified, all the (co)homology groups
%are in $\C$-coefficients.

\bigskip\noindent
{\bf Acknowledgments:} The authors thank Y. Ruan for stimulating discussions.
The third author also thanks Hong Kong UST for
its warm hospitality and support.
\section{\bf Preliminaries}

\subsection{Stable maps and Gromov-Witten invariants}
\label{subsect_gw}
\par
$\,$

Let $Y$ be a smooth projective variety.
A $k$-pointed {\it stable map} to $Y$ consists of
a complete nodal curve $C$ with $k$ distinct ordered smooth points
$p_1, \ldots, p_k$ and a morphism $\mu: C \to Y$ such that
the data $(\mu, C, p_1, \ldots, p_k)$ has only finitely many automorphisms.
In this case, the stable map is denoted by
$[\mu: (C; p_1, \ldots, p_k) \to Y]$.
For a fixed homology class $\beta \in H_2(Y; \mathbb Z)$,
let $\overline {\frak M}_{g, k}(Y, \beta)$ be the stack
parameterizing all the stable maps $[\mu: (C; p_1, \ldots, p_k) \to Y]$
such that $\mu_*[C] = \beta$ and the arithmetic genus of $C$ is $g$.
It is known \cite{F-P, LT1, LT2, B-F}
that $\overline {\frak M}_{g, k}(Y, \beta)$ is
a complete Deligne-Mumford stack with a projective moduli space.
Moreover, it has a virtual fundamental class
$[\overline {\frak M}_{g, k}(Y, \beta)]^{\text{\rm vir}} \in
A_{\frak d}(\overline {\frak M}_{g, k}(Y, \beta))$ where
\begin{eqnarray}  \label{expected-dim}
\frak d = -(K_Y \cdot \beta) + (\dim (Y) - 3)(1-g) + k
\end{eqnarray}
is the expected complex dimension of
$\overline {\frak M}_{g, k}(Y, \beta)$,
and $A_{\frak d}(\overline {\frak M}_{g, k}(Y, \beta))$
is the Chow group of $\frak d$-dimensional cycles in
the stack $\overline {\frak M}_{g, k}(Y, \beta)$.
The evaluation map
\begin{eqnarray}\label{evk}
ev_k\colon \overline {\frak M}_{g, k}(Y, \beta) \to Y^k
\end{eqnarray}
is defined by $ev_k([\mu: (C; p_1, \ldots, p_k) \to Y]) =
(\mu(p_1), \ldots, \mu(p_k))$.

The Gromov-Witten invariants are defined by using
the virtual fundamental class
$[\overline {\frak M}_{g, k}(Y, \beta)]^{\text{\rm vir}}$.
Recall that an element
$\alpha \in H^*(Y) {\buildrel\text{def}\over=}
\bigoplus_{j=0}^{2 \dim_{\mathbb C}(Y)} H^j(Y)$ is
{\it homogeneous} if $\alpha \in H^j(Y)$ for some $j$;
in this case, we take $|\alpha| = j$.
Let $\alpha_1, \ldots, \alpha_k \in H^*(Y)$
such that every $\alpha_i$ is homogeneous and
\begin{eqnarray}\label{homo-deg}
\sum_{i=1}^k |\alpha_i| = 2 {\frak d}.
\end{eqnarray}
Then, we have the $k$-point Gromov-Witten invariant defined by:
\begin{eqnarray}\label{def-GW}
\langle \alpha_1, \ldots, \alpha_k \rangle_{g, \beta} \,\,
= \int_{[\overline {\frak M}_{g, k}(Y, \beta)]^{\text{\rm vir}}}
ev_k^*(\alpha_1 \otimes \ldots \otimes \alpha_k).
\end{eqnarray}

Next, we summarize certain properties concerning the
virtual fundamental class. To begin with, we recall that
{\it the excess dimension} is the difference between the dimension of
$\overline {\frak M}_{g, k}(Y, \beta)$ and
the expected dimension $\frak d$ in (\ref{expected-dim}).
Let $T_Y$ stand for the tangent bundle of $Y$.
For $0 \le i < k$, we shall use
\begin{eqnarray}\label{k-to-i}
f_{k, i}: \overline {\frak M}_{g, k}(Y, \beta) \to
\overline {\frak M}_{g, i}(Y, \beta)
\end{eqnarray}
to stand for the forgetful map
obtained by forgetting the last $(k-i)$ marked points
and contracting all the unstable components.
It is known that $f_{k, i}$ is flat when $\beta \ne 0$ and $0 \le i < k$.
The following can be found in \cite{LT1, Beh, Get, C-K, LiJ}.

\begin{proposition}\label{virtual-prop}
Let $\beta \in H_2(Y; \mathbb Z)$ and $\beta \ne 0$.
Let $e$ be the excess dimension of
$\overline {\frak M}_{g, k}(Y, \beta)$, and $\frak M \subset
{\frak M}_{g, k}(Y, \beta)$ be a closed substack. Then,
\begin{enumerate}
\item[{\rm (i)}] $[\overline {\frak M}_{g, k}(Y, \beta)]^{\text{\rm vir}}
= (f_{k, 0})^*[\overline {\frak M}_{g, 0}(Y, \beta)]^{\text{\rm vir}}$;

%\item[{\rm (ii)}] $[\overline {\frak M}_{g, k}(Y, \beta)]^{\text{\rm vir}}
%= c_e(R^1(f_{k+1, k})_*(ev_{k+1})^*T_Y)$ if
%$R^1(f_{k+1, k})_*(ev_{k+1})^*T_Y$ is a rank-$e$ locally free sheaf
%over the moduli space $\overline {\frak M}_{g, k}(Y, \beta)$;

\item[{\rm (ii)}] $[\overline {\frak M}_{g, k}(Y, \beta)
]^{\text{\rm vir}}|_{\frak M}
= c_e((R^1(f_{k+1, k})_*(ev_{k+1})^*T_Y)|_{\frak M})$ if there exists
an open substack $\frak U$ of $\overline {\frak M}_{g, k}(Y, \beta)$
such that $\frak M \subset \frak U$ (i.e, $\frak U$ is
an open neighborhood of $\frak M$) and
$(R^1(f_{k+1, k})_*(ev_{k+1})^*T_Y)|_{\frak U}$ is
a rank-$e$ locally free sheaf over $\frak U$.
\end{enumerate}
\end{proposition}

We also need one formula for $g=0$ known as the composition law.
Let $\{ \Delta_a \}$ be a basis of $H^*(Y)$, and
$\{ \Delta^a \}$ be the basis of $H^*(Y)$ dual to $\{ \Delta_a \}$
with respect to the intersection pairing of $Y$.
Let $\alpha_1, \alpha_2, \alpha_3,
\alpha_4 \in H^*(Y)$ be classes of even degrees.
Then the combination of (3.3) and (3.6) in \cite{K-M} says that
\begin{eqnarray} \label{comp_law}
& &\langle \alpha_1 \alpha_2, \alpha_3, \alpha_4 \rangle_{0, \beta}
   +\langle \alpha_1, \alpha_2, \alpha_3 \alpha_4 \rangle_{0, \beta}
   \nonumber \\
& &\quad +\sum_{\beta_1 + \beta_2 = \beta, \, \beta_1, \beta_2 \ne 0} \,\,
   \sum_{a} \,\, \langle \alpha_1, \alpha_2, \Delta_a \rangle_{0, \beta_1}
   \, \langle \Delta^a, \alpha_3, \alpha_4 \rangle_{0, \beta_2}
   \nonumber \\
&=&\langle \alpha_1 \alpha_3, \alpha_2, \alpha_4 \rangle_{0, \beta}
   +\langle \alpha_1, \alpha_3, \alpha_2 \alpha_4 \rangle_{0, \beta}
   \nonumber \\
& &\quad +\sum_{\beta_1 + \beta_2 = \beta, \, \beta_1, \beta_2 \ne 0} \,\,
   \sum_{a} \,\, \langle \alpha_1, \alpha_3, \Delta_a \rangle_{0, \beta_1}
   \, \langle \Delta^a, \alpha_2, \alpha_4 \rangle_{0, \beta_2}.
\end{eqnarray}
%the matrix $\{ g^{ab} \}$ be the inverse of the intersection matrix of
%the basis $\{ \Delta_a \}$. Let $\alpha_1, \alpha_2, \alpha_3,
%\alpha_4 \in H^*(Y)$ be classes of even degrees.
%Then the combination of (3.3) and (3.6) in \cite{K-M} says that
%\begin{eqnarray} \label{comp_law}
%& &\langle \alpha_1 \alpha_2, \alpha_3, \alpha_4 \rangle_{0, \beta}
%   +\langle \alpha_1, \alpha_2, \alpha_3 \alpha_4 \rangle_{0, \beta}
%   \nonumber \\
%& &\quad +\sum_{\beta_1 + \beta_2 = \beta, \, \beta_1, \beta_2 \ne 0} \,\,
%   \sum_{a, b} \,\, \langle \alpha_1, \alpha_2, \Delta_a \rangle_{0, \beta_1}
%   \,\, g^{ab} \,\,
%   \langle \Delta_b, \alpha_3, \alpha_4 \rangle_{0, \beta_2}
%   \nonumber \\
%&=&\langle \alpha_1 \alpha_3, \alpha_2, \alpha_4 \rangle_{0, \beta}
%   +\langle \alpha_1, \alpha_3, \alpha_2 \alpha_4 \rangle_{0, \beta}
%   \nonumber \\
%& &\quad +\sum_{\beta_1 + \beta_2 = \beta, \, \beta_1, \beta_2 \ne 0} \,\,
%   \sum_{a, b} \,\, \langle \alpha_1, \alpha_3, \Delta_a \rangle_{0, \beta_1}
%   \,\, g^{ab} \,\,
%   \langle \Delta_b, \alpha_2, \alpha_4 \rangle_{0, \beta_2}.
%\end{eqnarray}

%
%
%
%
%
%
%
\subsection{Basic facts about the Hilbert scheme of points on a surface}
\label{subsect_fact}
\par
$\,$

Let $X$ be a simply connected smooth projective surface,
and $\Xn$ be the Hilbert scheme of points in $X$.
An element in $\Xn$ is represented by a length-$n$ $0$-dimensional
closed subscheme $\xi$ of $X$. For $\xi \in \Xn$,
let $I_{\xi}$  be the corresponding sheaf of ideals.
In $\Xn\times X$, we have the universal codimension-$2$ subscheme:
\begin{eqnarray}\label{univ-subscheme}
\mathcal Z_n=\{(\xi, x) \subset \Xn\times X \, |\,
x\in \Supp{(\xi)}\}\subset \Xn\times X.
\end{eqnarray}
%In $X^{[n-1]}\times \Xn$, we have the $2n$-dimensional smooth
%incidence subscheme:
%\begin{eqnarray}\label{incidence}
%X^{[n-1, n]}=\{ (\xi,\eta)\in X^{[n-1]}\times \Xn \,|\,
%I_{\xi}\supset I_{\eta}\}.
%\end{eqnarray}
Let $X^{(n)}$ be the $n$-th symmetric product of $X$.
We have the Hilbert-Chow map:
\begin{eqnarray}  \label{H_C}
\rho: X^{[n]} \to X^{(n)}.
\end{eqnarray}

For a subset $Y \subset X$, we define the subset $M_n(Y)$ in
the Hilbert scheme $\Xn$:
\begin{eqnarray}\label{subset}
M_n(Y) = \{ \xi \in \Xn| \Supp(\xi) \text{ is a point in } Y \}
\subset \Xn.
\end{eqnarray}
In particular, for a fixed point $x \in X$,
$M_n(x)$ is just the punctual Hilbert scheme of points on $X$ at $x$.
It is known that the punctual Hilbert schemes $M_n(x)$ are
isomorphic for all the surfaces $X$ and all the points $x \in X$.

%The definitions and properties of the maps listed below can
%be found in \cite{ES2}.
%
%\begin{notation}\label{notations}
%There exist various morphisms:
%\begin{eqnarray*}
%&f_n&\colon \,\, X^{[n-1, n]} \to X^{[n-1]}
% \text{ with } f_n(\xi, \eta)=\xi.\\
%&g_n&\colon \,\, X^{[n-1, n]}\to \Xn
% \text{ with } g_n(\xi, \eta)=\eta.\\
%&\psi_n&\colon \,\, X^{[n-1, n]}\to \mathcal Z_n
% \text{ with } \psi_n(\xi, \eta) =
%      \big(\eta, \Supp(I_{\xi}/I_{\eta})\big).\\
%&q&\colon \,\, X^{[n-1, n]}\to X
% \text{ with } q(\xi, \eta)=\Supp(I_{\xi}/I_{\eta}).
%\end{eqnarray*}
%\end{notation}
%
%\begin{convention}
%Let $V$ be an $n$-dimensional vector space. We use
%$\Bbb P(V)$ to denote the set of $1$-dimensional quotients of
%the vector space $V$.
%\end{convention}
%
%\begin{theorem}\label{ES'sthm} (see \cite{ES2}) Adopt the above notations.
%\begin{enumerate}
%\item[{\rm (i)}] The morphism $\psi_{n}\colon\, X^{[n-1, n]}\rightarrow
%\mathcal Z_{n}$ is canonically isomorphic to the projectification
%$\mathbb P(\omega_{\mathcal Z_{n}}) \to \mathcal Z_{n}$ where
%$\omega_{\mathcal Z_{n}}$ is the dualizing sheaf of $\mathcal Z_{n}$;
%
%\item[{\rm (ii)}] The morphism $(f_n, q) \colon\, X^{[n-1, n]}\rightarrow
%X^{[n-1]} \times X$ is canonically isomorphic to the blowing-up of
%$X^{[n-1]} \times X$ along $\mathcal Z_{n-1}$. The exceptional locus is
%\begin{eqnarray}\label{ex-locus}
%E_n=\{ (\xi, \eta)\in X^{[n-1, n]}\, |\,
%\Supp(\xi)=\Supp(\eta) \text{ and } \xi\subset \eta\}.
%\end{eqnarray}
%\end{enumerate}
%\end{theorem}

Let $\xi\in X^{[n-k]}$ and $\eta\in X^{[k]}$.
If $\Supp(\xi)\cap \Supp(\eta)=\emptyset$, then we use $\xi+\eta$ to
represent the closed subscheme $\xi\cup\eta$ in $\Xn$. Similarly,
given a subvariety $Y$ of $X^{[n-k]}$ and a point $\eta\in X^{[k]}$
such that $\left (\bigcup\limits_{\xi\in Y}\Supp(\xi)\right )\cap
\Supp(\eta) =\emptyset$, we use $Y+\eta$ to represent
the subvariety in $\Xn$ consisting of all the points
$\xi + \eta$ with $\xi \in Y$.

Next, we review some results on homology groups of the Hilbert scheme
$\Xn$ due to G\"ottsche \cite{Go1}, Grojnowski \cite{Gro},
and Nakajima \cite{Nak}. Their results say that the space
$\mathbb H \buildrel\text{def}\over=
\bigoplus\limits^{\infty}_{n=0}\bigoplus\limits^{4n}_{k=0}
H_k(\Xn)$ is an irreducible highest weight representation of
the Heisenberg algebra generated by
$\frak a_{-n}(a), n \in \mathbb Z, a \in H_*(X)
\buildrel\text{def}\over= \bigoplus\limits^{4}_{k=0} H_k(X)$.
Moreover, $\vac \buildrel\text{def}\over= 1\in H_0(X^{[0]};
\mathbb C)=\mathbb C$ is a highest weight vector. It follows that
the space $\mathbb H$ is a linear span of elements of the form
$\frak a_{-n_1}(a_1) \ldots\frak a_{-n_k}(a_k) \vac$
where $k \ge 0$, $n_1, \ldots, n_k > 0$, and
$a_1, \ldots, a_k \in H_*(X)$. The geometric interpretation of
$\frak a_{-n_1}(a_1) \ldots\frak a_{-n_k}(a_k) \vac$ for homogeneous
classes $a_1, \ldots, a_k \in H_*(X)$ can be understood as follows.
For $i = 1, \ldots, k$, let $a_i \in H_{|a_i|}(X)$
be represented by a cycle $X_i$ such that
$X_1, \ldots, X_k$ are in general position. Then,
\begin{eqnarray}\label{bidegree}
\frak a_{-n_1}(a_1) \ldots\frak a_{-n_k}(a_k)\vac \in H_m(\Xn)
\end{eqnarray}
where $n=\sum\limits_{i=1}^k n_i$ and
$m = \sum\limits_{i=1}^k(2n_i-2+|a_i|)$.
Up to a scalar, $\frak a_{-n_1}(a_1)
\ldots\frak a_{-n_k}(a_k)\vac$ is represented by the closure of
the real-$\sum\limits_{i=1}^k(2n_i-2+|a_i|)$-dimensional subset:
\begin{eqnarray} \label{geom-rep}
\{ \xi_1+\ldots + \xi_k \in \Xn|
\xi_i \in M_{n_i}(X_i), \Supp(\xi_i) \cap \Supp(\xi_j) =
\emptyset \text{ for } i \ne j \}
\end{eqnarray}
where $M_{n_i}(X_i)$ is the subset of
$X^{[n_i]}$ defined by (\ref{subset}).

\begin{definition}\label{cohomology-classes}
Let $x \in X$, and $C$ be a real-$2$-dimensional submanifolds of $X$.
Then, we define $\beta_n = \frak a_{-2}(x) \frak a_{-1}(x)^{n-2}\vac$,
$\beta_C = \frak a_{-1}(C)\frak a_{-1}(x)^{n-1}\vac$, and
\begin{eqnarray*}
B_n = {1 \over (n-2)!} \frak a_{-2}(X) \frak a_{-1}(X)^{n-2}\vac,
%\quad{\text{\rm and}}
\quad
D_C = {1 \over (n-1)!} \frak a_{-1}(C) \frak a_{-1}(X)^{n-1}\vac.
\end{eqnarray*}
\end{definition}

\begin{lemma}\label{lemmacohoclass}
Let $x$ and $\ell$ be a point and a line in $X = \mathbb P^2$
respectively. Then,
\begin{enumerate}
\item[{\rm (i)}] a basis of $H_2(X^{[3]}; \mathbb Z)$ consists of
$\beta_3$ and $\beta_{\ell}$;

\item[{\rm (ii)}] a basis of $H_4(X^{[3]})$ consists of
the five homology classes
$\frak a_{-1}(X)\frak a_{-1}(x)^{2} \vac$,
$\frak a_{-2}(\ell)\frak a_{-1}(x)\vac$,
$\frak a_{-1}(\ell)^2 \frak a_{-1}(x)\vac$,
$\frak a_{-1}(\ell)\frak a_{-2}(x)\vac$,
and $\frak a_{-3}(x)\vac$;

\item[{\rm (iii)}] a basis of $H_6(X^{[3]})$ consists of
the classes
$\frak a_{-2}(X)\frak a_{-1}(x) \vac$,
$\frak a_{-1}(X)\frak a_{-2}(x) \vac$,
$\frak a_{-1}(X)\frak a_{-1}(\ell)\frak a_{-1}(x)\vac$,
$\frak a_{-3}(\ell)\vac$, $\frak a_{-2}(\ell)\frak a_{-1}(\ell)\vac$,
and $\frak a_{-1}(\ell)^3\vac$;

\item[{\rm (iv)}] a basis of $H_8(X^{[3]})$ consists of
the five classes $\frak a_{-3}(X)\vac$,
$\frak a_{-2}(X)\frak a_{-1}(\ell)\vac$,
$\frak a_{-1}(X)\frak a_{-2}(\ell)\vac$,
$\frak a_{-1}(X)\frak a_{-1}(\ell)^2\vac$,
and $\frak a_{-1}(X)^{2}\frak a_{-1}(x) \vac$;

\item[{\rm (v)}] a basis of $H_{10}(X^{[3]}; \mathbb Z)$ consists of
the divisors $B_3$ and $D_{\ell}$.
\end{enumerate}
\end{lemma}
\begin{proof}
The proof of (i) and (v) was contained in the proof of
the Theorem 4.1 in \cite{LQZ}, while the rest statements follow by
exploiting (\ref{bidegree}).
\end{proof}

\begin{definition} \label{def_B_i}
For $X = \mathbb P^2$ and $i = 2, 4, 6, 8$ and $10$,
let $\frak B_i$ stand for the linear basis of the homology group
$H_i(X^{[3]})$ given in Lemma~\ref{lemmacohoclass}.
\end{definition}

Fix $p \in X^{[3]}$. Then a basis $\{ \Delta_a \}$ of
$H^*(X^{[3]})$ is given by the Poincar\'e duals of
\begin{eqnarray} \label{basis}
[p], \,\, \frak B_i \, (i = 2, 4, 6, 8, 10), \,\, [X^{[3]}]
\end{eqnarray}
where $[p] = \frak a_{-1}(x)^3 \vac \in H_0(X^{[3]})$ and
$[X^{[3]}] = 1/6 \,\, \frak a_{-1}(X)^3 \vac \in H_{12}(X^{[3]})$ are
the homology classes corresponding to $p$ and $X^{[3]}$ respectively.

The following is the main result proved in \cite{L-Q}.

\begin{lemma} \label{1pt}
Let $d\ge 1$, and $x$ and $\ell$ be a point and a line
in $X = \mathbb P^2$ respectively.

{\rm (i)} If $\alpha$ stands for the Poincar\'e duals of
the homology classes
$\frak a_{-1}(X)\frak a_{-1}(x)^{2} \vac$,
$\frak a_{-1}(\ell)^2 \frak a_{-1}(x)\vac$,
$\frak a_{-1}(\ell)\frak a_{-2}(x)\vac$,
and $\frak a_{-3}(x)\vac$,
then $\langle \alpha\rangle_{0, d\beta_n}=0$.

{\rm (ii)} If $\alpha$ is the Poincar\'e dual of
$\frak a_{-2}(\ell)\frak a_{-1}(x)\vac$, then
$\langle \alpha\rangle_{0, d\beta_n}=2(K_X\cdot \ell)/d^2$.
\end{lemma}

\subsection{Curves from the punctual Hilbert scheme}
\label{subsect_punctual}
\par
$\,$

\begin{lemma}     \label{sigma_n}
Fix $n \ge 2$.
Let $\Hilb^n(\mathbb C^2, 0)$ be the punctual Hilbert scheme of
points on $\mathbb C^2$ at the origin,
and $u, v$ be the coordinates of $\mathbb C^2$. Then,
$H_2(\Hilb^n(\mathbb C^2, 0); \mathbb Z) \cong \mathbb Z$.
Moreover, a generator of $H_2(\Hilb^n(\mathbb C^2,0); \mathbb Z)$
is given by
\begin{eqnarray}  \label{sigma_n:def}
\sigma_n=\{(\lambda u + \mu v^{n-1}, u^2, uv, v^n)\,| \,
\lambda, \mu \in \mathbb C \text{ with } |\lambda|+|\mu| \neq 0 \}.
\end{eqnarray}
\end{lemma}
\begin{proof}
The first statement was proved in \cite{E-S}.
To prove the second statement, following \cite{E-S},
take a $\mathbb C^*$-action on $\mathbb C^2$ given by
$t\cdot(u, v)=(t^{-\alpha}u, t^{-\beta}v)$ with $\beta\gg \alpha$.
For $\xi \in \Hilb^n(\C^2; 0)$,
we use the ideal $I_{\xi} \subset \C[u,v]$ to represent $\xi$.
Then the $\C^*$-invariant ideal in
$\C[u, v]$ corresponding to a generator $\sigma_n$ of
$H_2(\Hilb^n(\C^2, 0); \Z)$ is $(v^{n-1}, uv, u^2)$.
Therefore $\sigma_n$ is the closure of the cell
\begin{eqnarray*}
& &\{I\in \C[u, v]\,|\,\ell({\C[u, v]/I})=n,
  \quad\lim_{t\to 0}(t\cdot I) =(v^{n-1}, uv,u^2)\} \\
&=& \{(v^{n-1}+au, uv, u^2) \,| \, a\in \C\} \cong \C.
\end{eqnarray*}
Finally, notice that if $a \ne 0$, then $(v^{n-1}+au, uv, u^2)
=(v^{n-1}+au, v^n)$. So letting $a \to \infty$, we see that
the ideal $(u, v^n)$ is also contained in $\sigma_n$. Thus,
\begin{eqnarray*}
\sigma_n =\{(v^{n-1}+au, uv, u^2) \,| \, a\in \C\}\cup\{(u, v^n)\}
\end{eqnarray*}
which is the same as $\{(\lambda u + \mu v^{n-1}, u^2, uv, v^n)
\,| \, \lambda, \mu \in \C \text{ with } |\lambda|+|\mu| \neq 0 \}$.
\end{proof}

Let $R = \mathcal O_{\C^2, 0}$ be the local ring of $\C^2$ at
the origin, and $\mathfrak m=(u, v)$ be the maximal ideal of $R$.
Let $\eta \in \Hilb^n(\C^2, 0)$.
It is known that there exists an embedding
\begin{eqnarray*}
\tau: \Hilb^n(\C^2, 0) \to \Grass(R/\mathfrak m^n, n)
\end{eqnarray*}
where $R/\mathfrak m^n$ is considered as a $\C$-vector space
of dimension $n+1\choose 2$, and $\tau$ maps
an element $\eta\in \Hilb^n(\C^2, 0)$ to the $n$-dimensional quotient
of $R/\mathfrak m^n$ in the exact sequence
\begin{eqnarray*}
0 \to I_{\eta, 0}/\mathfrak m^n \to R/\mathfrak m^n \to
R/I_{\eta, 0}=\mathcal O_{\eta, 0} \to 0.
\end{eqnarray*}
Let $\mathfrak p: \mathbb G \to \mathbb P^{N-1}$ be the Pl\"ucker
embedding where $N = {n+1\choose 2}\left ({n+1\choose 2}-n \right )$.

\begin{lemma} \label{sigma_n:homology}
Identify $M_n(x)$ with $\Hilb^n(\C^2, 0)$, and regard
$\sigma_n$ as a curve in $M_n(x) \subset \Xn$.
Then as a curve in $\Xn$, $\sigma_n$ is homologous to $\beta_n$.
\end{lemma}
\begin{proof}
According to the results in Sect. 3 of \cite{LQZ}, it suffices to show
that the image $(\mathfrak p \circ \tau)(\sigma_n)$ is a line.
Fix a basis for the $\mathbb C$-vector space $R/\mathfrak m^n$:
\begin{eqnarray*}
\overline 1, \overline u, \overline u^2, \overline u \, \overline v,
\overline u^3, \overline u^2 \, \overline v,
\overline u  \, \overline v^2, \ldots, \overline u^{n-1},
\overline u^{n-2} \overline v, \ldots, \overline u  \, \overline v^{n-2},
\overline v, \ldots, \overline v^{n-1}.
\end{eqnarray*}
Note the special ordering of this basis.
Recall from (\ref{sigma_n:def}) that
for any $\eta \in \sigma_n \subset \Hilb^n(\C^2, 0)$,
$I_{\eta, 0}=(\lambda u + \mu v^{n-1}, u^2, uv, v^n)$
for some $\lambda, \mu\in \mathbb C$ with $|\lambda|+|\mu| \neq 0$.
So a basis for the subspace $I_{\eta, 0}/\mathfrak m^n \subset R/\mathfrak m^n$
can be chosen as
\begin{eqnarray*}
\lambda \overline u+\mu \overline v^{n-1}, \overline u^2,
\overline u \, \overline v, \overline u^3, \overline u^2 \overline v,
\overline u \, \overline v^2, \ldots, \overline u^{n-1},
\overline u^{n-2} \overline v, \ldots,
\overline u \, \overline v^{n-2},
\end{eqnarray*}
and the matrix representation of $I_{\eta, 0}/\frak m^n$ is given by
the ${n\choose 2}\times{ n+1\choose 2}$-matrix:
\begin{eqnarray}
 \left[
  \begin{array}{ccccccccc}
0&\lambda&0&\ldots&0&0&\ldots&0&\mu\\
0&0&1&\ldots&0&0&\ldots&0&0\\
\vdots&\vdots&\vdots&\vdots&\vdots&\vdots&\ldots&\vdots&\vdots\\
0&0&0&\ldots&1&0&\ldots&0&0
  \end{array}
 \right].
\end{eqnarray}
%\begin{eqnarray}
% \left[
%  \begin{array}{cc}
%   f_0 (D)  &   f_+  (D)   \\
%   f_- (D)  &   f_1  (D)
%  \end{array}
% \right],
%\end{eqnarray}
%You can change {cc} to {ccc} to get 3X3 matrix, etc.
%BTW, "c" here stands for "center"; one can change
%"cc" to "ll" or "rr" to move the entries of a
%matrix to the left and right.
Thus, $(\mathfrak p\circ \tau)(\eta)=[0, \ldots, 0, \lambda, 0, \ldots,
0, \mu, 0, \ldots, 0]$ where the positions of $\lambda$ and $\mu$
are independent of $\eta\in \sigma_n$.
So the image $(\frak p\circ \tau)(\sigma_n)$ is a line.
\end{proof}

Note that the flat limits of the elements $(\lambda u + v, v^n)$,
$\lambda \in \C^*$ in $\Hilb^n(\C^2, 0)$ as $\lambda \to 0$
and $\lambda \to \infty$ are equal to $(v, u^n)$ and $(u, v^n)$ respectively.
So in the punctual Hilbert scheme $\Hilb^n(\C^2, 0)$,
we have the projective curve:
\begin{eqnarray}   \label{tilde_sigma_n:def}
\tilde \sigma_n=\{(\lambda u + v, v^n)\,| \, \lambda \in \C^* \}
\cup \{ (v, u^n), \,\,  (u, v^n)\}.
\end{eqnarray}

\begin{lemma} \label{tilde_sigma_n:homology}
As a curve in $\Xn$, $\tilde \sigma_n$ is homologous to
${n\choose 2} \sigma_n$.
\end{lemma}
\begin{proof}
It suffices to show that $\tilde \sigma_n \sim {n\choose 2} \sigma_n$
in $H_2(\Hilb^n(\C^2, 0); \mathbb Z)$. By (\ref{tilde_sigma_n:def}),
if $\eta \in \tilde \sigma_n - \{ (v, u^n), \,\,  (u, v^n) \}$,
then a basis for the subspace
$I_{\eta, 0}/\mathfrak m^n \subset R/\mathfrak m^n$ is
\begin{eqnarray*}
&\lambda \overline u+ \overline v,
\lambda \overline u^2+ \overline u \, \overline v,
\lambda \overline u \, \overline v + \overline v^2, \ldots,& \\
&\lambda \overline u^{n-1} + \overline u^{n-2} \overline v,
\lambda \overline u^{n-2} \overline v + \overline u^{n-3} \overline v^2,
\ldots, \lambda \overline u \, \overline v^{n-2} + \overline v^{n-1}.&
\end{eqnarray*}
As in the proof of Lemma~\ref{sigma_n:homology}, we see that the degree
of $(\frak p\circ \tau)(\tilde \sigma_n - \{ (v, u^n), \,\,  (u, v^n) \})$
is ${n\choose 2}$.
So $(\frak p\circ \tau)(\tilde \sigma_n)$ has degree ${n\choose 2}$.
By Lemma~\ref{sigma_n}, there exists an integer $d$ such that
$\tilde \sigma_n \sim d \sigma_n$ in $H_2(\Hilb^n(\C^2, 0); \mathbb Z)$.
Since $(\frak p\circ \tau)(\sigma_n)$ is a line, $d = {n\choose 2}$.
\end{proof}
\section{\bf $3$-point genus-$0$ Gromov-Witten invariants of
$(\PP^2)^{[3]}$}
\label{sect_GW}

Let $X= \PP^2$ and $d \ge 1$. For simplicity, we shall use
$\langle \alpha_1, \ldots, \alpha_k \rangle_{0, d}$ to stand for
$\langle \alpha_1, \ldots, \alpha_k \rangle_{0, d\beta_3}$.
Our goal is to compute the $3$-point Gromov-Witten
invariants $\langle \alpha_1, \alpha_2, \alpha_3 \rangle_{0, d}$
of $X^{[3]}$. Recall from Lemma~\ref{1pt} that
the $1$-point Gromov-Witten invariants
$\langle \alpha_1 \rangle_{0, d}$ of $X^{[3]}$ have been calculated.
Since the expected complex dimension of
the stack $\overline {\frak M}_{0, 3}(X^{[3]}, d\beta_3)$ is $6$,
it remains to compute the $2$-point Gromov-Witten invariants
$\langle \pd(A_1), \pd(A_2) \rangle_{0, d}$
when $A_1$ runs over the basis $\frak B_6$ of $H_6(X^{[3]})$ in
Lemma~\ref{lemmacohoclass}~(iii) and $A_2$ runs over the basis
$\frak B_8$ of $H_8(X^{[3]})$ in Lemma~\ref{lemmacohoclass}~(iv),
and $\langle \pd(A_1), \pd(A_2), \pd(A_3) \rangle_{0, d}$
when $A_1, A_2, A_3$ run over the basis $\frak B_8$.
\subsection{$\langle \pd(A_1), \pd(A_2) \rangle_{0, d}$
with $A_1 \in \frak B_6$ and $A_2 \in \frak B_8$}
\label{subsect_2}
\par
$\,$

\begin{lemma} \label{2_zero}
The $2$-point Gromov-Witten invariants
$\langle \pd(A_1), \pd(A_2) \rangle_{0, d}$ are equal to zero
for the following pairs of $(A_1, A_2) \in
\frak B_6 \times \frak B_8$:
\begin{eqnarray*}
&(\frak a_{-2}(X)\frak a_{-1}(x) \vac,
  \frak a_{-2}(X)\frak a_{-1}(\ell)\vac),\,\,
  (\frak a_{-1}(X)\frak a_{-2}(x) \vac,
  \frak a_{-1}(X)\frak a_{-2}(\ell)\vac),& \\
&(\frak a_{-1}(X)\frak a_{-2}(x) \vac,
  \frak a_{-1}(X)\frak a_{-1}(\ell)^2\vac), \,\,
  (\frak a_{-1}(X)\frak a_{-1}(\ell)\frak a_{-1}(x)\vac,
  \frak a_{-3}(X)),& \\
&(\frak a_{-1}(X)\frak a_{-1}(\ell)\frak a_{-1}(x)\vac,
  \frak a_{-1}(X)\frak a_{-1}(\ell)^2\vac), &\\
&(\frak a_{-1}(X)\frak a_{-1}(\ell)\frak a_{-1}(x)\vac,
  \frak a_{-1}(X)^{2}\frak a_{-1}(x) \vac),\,\,
 (\frak a_{-3}(\ell)\vac,
  \frak a_{-1}(X)\frak a_{-1}(\ell)^2\vac),&\\
&(\frak a_{-3}(\ell)\vac,
  \frak a_{-1}(X)^{2}\frak a_{-1}(x) \vac), \,\,
 (\frak a_{-2}(\ell)\frak a_{-1}(\ell)\vac,
  \frak a_{-1}(X)^{2}\frak a_{-1}(x) \vac), \\
&(\frak a_{-1}(\ell)^3\vac, \frak a_{-3}(X)),\,\,
  (\frak a_{-1}(\ell)^3\vac,
  \frak a_{-1}(X)\frak a_{-2}(\ell)\vac),&\\
&(\frak a_{-1}(\ell)^3\vac,
  \frak a_{-1}(X)\frak a_{-1}(\ell)^2\vac), \,\,
  (\frak a_{-1}(\ell)^3\vac,
  \frak a_{-1}(X)^{2}\frak a_{-1}(x) \vac).&
\end{eqnarray*}
\end{lemma}
\begin{proof}
These follow from similar geometric arguments. For instance, let us
show that $\langle \pd(A_1), \pd(A_2) \rangle_{0, d} = 0$
when $A_1=\frak a_{-1}(\ell)^3\vac$ and
$A_2=\frak a_{-1}(X)\frak a_{-1}(\ell)^2\vac$.

Choose five lines $\ell_1, \ldots, \ell_5 \subset X = \PP^2$
in general position. By (\ref{geom-rep}), we see that up to a scalar,
$A_1$ is represented by the closure of the subset
\begin{eqnarray}  \label{2_zero.1}
\{ x_1+x_2+x_3|\,\, x_1, x_2, x_3 \,\, \text{are distinct and }
x_i \in \ell_i \,\, \text{for each } i\}.
\end{eqnarray}
Similarly, $A_2$ is represented by the closure of the subset
\begin{eqnarray}  \label{2_zero.2}
\{ x+x_4+x_5|\,\, x, x_4, x_5 \,\, \text{are distinct and }
x_i \in \ell_i \,\, \text{for each } i\}.
\end{eqnarray}

Let $\frak M$ be the substack of
$\overline {\frak M}_{0, 2}(X^{[3]}, d\beta_3)$
parametrizing all the stable maps $[\mu: (C; p_1, p_2) \to X^{[3]}]$
with $\mu(p_1) \in A_1$ and $\mu(p_2) \in A_2$.
We claim that $\frak M = \emptyset$. Indeed, assume
$[\mu: (C; p_1, p_2) \to X^{[3]}]$ is an object of $\frak M$. On one hand,
by (\ref{2_zero.1}), $\rho(\mu(C)) = 2(\ell_i \cap \ell_j) + x_k$
where $\rho$ is the Hilbert-Chow map~(\ref{H_C}),
$\{i, j, k\}$ is a permutation of $\{1, 2, 3\}$, and $x_k \in \ell_k$.
On the other hand, by (\ref{2_zero.2}), we obtain
\begin{eqnarray*}
\rho(\mu(C)) = 2(\ell_4 \cap \ell_5) + x
\end{eqnarray*}
for some $x \in X$, or $\rho(\mu(C)) = 2x_i + x_j$ where $\{i, j\}$ is
a permutation of $\{4, 5\}$, $x_i \in \ell_i$, and $x_j \in \ell_j$.
Since the lines $\ell_1, \ldots, \ell_5 \subset X = \PP^2$ are
in general position, such $\rho(\mu(C))$ does not exist.
So $\frak M = \emptyset$. Hence
$\langle \pd(A_1), \pd(A_2) \rangle_{0, d} = 0$.
\end{proof}

\begin{lemma} \label{2_zero_K3}
The $2$-point Gromov-Witten invariants
$\langle \pd(A_1), \pd(A_2) \rangle_{0, d}$ are equal to zero for
the following pairs of $(A_1, A_2) \in \frak B_6 \times \frak B_8$:
\begin{eqnarray*}
&(\frak a_{-2}(X)\frak a_{-1}(x) \vac,
  \frak a_{-3}(X)\vac),\,\,
  (\frak a_{-2}(X)\frak a_{-1}(x) \vac,
  \frak a_{-1}(X)\frak a_{-1}(\ell)^2\vac),\\
&(\frak a_{-2}(X)\frak a_{-1}(x) \vac,
  \frak a_{-1}(X)^{2}\frak a_{-1}(x) \vac), \,\,
  (\frak a_{-1}(X)\frak a_{-2}(x) \vac, \frak a_{-3}(X)\vac),& \\
&(\frak a_{-1}(X)\frak a_{-2}(x) \vac,
  \frak a_{-1}(X)^{2}\frak a_{-1}(x) \vac), \,\,
  (\frak a_{-1}(X)\frak a_{-1}(\ell)\frak a_{-1}(x)\vac,
  \frak a_{-2}(X)\frak a_{-1}(\ell)\vac),& \\
&(\frak a_{-1}(X)\frak a_{-1}(\ell)\frak a_{-1}(x)\vac,
  \frak a_{-1}(X)\frak a_{-2}(\ell)\vac), \,\,
  (\frak a_{-3}(\ell)\vac,\frak a_{-2}(X)\frak a_{-1}(\ell)\vac),&  \\
&(\frak a_{-3}(\ell)\vac,\frak a_{-1}(X)\frak a_{-2}(\ell)\vac),\,\,
  (\frak a_{-2}(\ell)\frak a_{-1}(\ell)\vac, \frak a_{-3}(X)\vac),& \\
&(\frak a_{-2}(\ell)\frak a_{-1}(\ell)\vac,
  \frak a_{-1}(X)\frak a_{-1}(\ell)^2\vac),\,\,
  (\frak a_{-1}(\ell)^3\vac, \frak a_{-2}(X)\frak a_{-1}(\ell)\vac).&
\end{eqnarray*}
\end{lemma}
\begin{proof}
These invariants are equal to certain genus-$0$ Gromov-Witten invariants
of a K3 surface. So our lemma follows from the fact that all the
genus-$0$ Gromov-Witten invariants of a K3 surface are equal to zero.
For instance, let us show that
$\langle \pd(A_1), \pd(A_2) \rangle_{0, d} = 0$
when $A_1=\frak a_{-2}(X)\frak a_{-1}(x) \vac$ and
$A_2=\frak a_{-3}(X)\vac$.

Fix $x \in X$, and a small analytic open subset $U$ of $X$
such that $x \in U$.
We may assume that $U$ is independent of $X$.
Note that for a stable map $[\mu: (C; p_1, p_2) \to X^{[3]}]
\in \overline{\frak M}_{0, 2}(X^{[3]}, d\beta_3)$,
either $\mu(C) \subset U^{[3]}$ or $\mu(C) \cap U^{[3]} = \emptyset$.
So the analytic open substack $\frak U \subset
\overline{\frak M}_{0, 2}(X^{[3]}, d\beta_3)$ parametrizing
all stable maps $[\mu\colon (C; p_1, p_2)\to X^{[3]}]$ with
$\mu(C) \subset U^{[3]}$ depends only on $U$,
and is independent of $X$.

Let $\frak M$ be the substack of
$\overline {\frak M}_{0, 2}(X^{[3]}, d\beta_3)$
parametrizing all the stable maps $[\mu: (C; p_1, p_2) \to X^{[3]}]$
such that $\mu(p_1) \in A_1$ and $\mu(p_2) \in A_2$.
Note from the descriptions of $A_1$ and $A_2$ that
if $[\mu: (C; p_1, p_2) \to X^{[3]}] \in \frak M$,
then $\mu(C) \subset M_3(x) \subset U^{[3]}$.
So $\frak M \subset \frak U$. In fact, $\frak M$ parametrizes
all the stable maps $[\mu: (C; p_1, p_2) \to X^{[3]}] \in \frak U$
with $\mu(C) \subset M_3(x) \subset U^{[3]}$.
So $\frak M$ is also independent of $X$.

In summary, we showed that $\frak M \subset \frak U$ where $\frak U$
is analytic open in $\overline{\frak M}_{0, 2}(X^{[3]}, d\beta_3)$,
and $\frak M$ and $\frak U$ are independent of $X$.
It follows from the constructions of the virtual fundamental class
(see \cite{LT2, LT3, Ru1}) that the restriction
$[\overline {\frak M}_{0, 2}(X^{[3]}, d\beta_3)]^{\text{\rm vir}}|_{\frak M}$
is independent of the smooth surface $X$. So we have
$\langle \pd(A_1), \pd(A_2) \rangle_{0, d} =
\langle \pd(A_1'), \pd(A_2') \rangle_{0, d}$
where $A_1'=\frak a_{-2}(X')\frak a_{-1}(x') \vac$,
$A_2'=\frak a_{-3}(X')\vac$, $x' \in X'$, and $X'$ is a K3 surface.
Therefore, we conclude that
$\langle \pd(A_1), \pd(A_2) \rangle_{0, d} =0$.
\end{proof}

To compute other $2$-point invariants
$\langle \pd(A_1), \pd(A_2) \rangle_{0, d}$,
we recall from \cite{L-Q} some results concerning obstruction bundles
and virtual fundamental classes.
Fix $n \ge 2$. Let $B_*=\{\xi\in \Xn\, |\, |\Supp(\xi)|=n-1\}$ and
$X^{(n)}_{s*} = \rho(B_*)$ where $\rho$ is the Hilbert-Chow map.
Let $j_2: X^{(n)}_{s*} \to X$ be the morphism defined by
sending $2x+x_3+\ldots+x_n$ to $x$.
For $k \ge 0$, let $\frak U_k$ be the open substack
of $\overline{\frak M}_{0, k}(\Xn, d\beta_n)$ parametrizing
stable maps $[\mu\colon (C; p_1, \ldots, p_k)\to \Xn]$
such that $\mu(C)\subset B_*$. For $k \ge 1$, note that
$\frak U_k=f^{-1}_{k, 0}(\frak U_0)$.
Put ${\tilde {ev}_k} = ev_k|_{\frak U_k}$ and
${\tilde {f}_{k, 0}} = f_{k, 0}|_{\frak U_k}$.
Then we can regard ${\tilde {ev}_k}$ and ${\tilde {f}_{k, 0}}$
as morphisms from $\frak U_k$ to $(B_*)^k$ and $\frak U_0$ respectively.
In addition, there exist morphisms $\phi$ and $j_1$ forming
a commutative diagram:
\begin{eqnarray}\label{com-diagram2}
\begin{matrix}
{\frak U_1}&{\buildrel{\rev}\over \rightarrow }  &B_*
  &{\buildrel {j_1}\over  \cong}&\mathbb P(j_2^*T_X^*)&&\\
\quad \downarrow^{\rf}&&\downarrow^{\rho}&&\downarrow^{\pi}&&\\
{\frak U_0}&{\buildrel{\phi}\over \rightarrow }&\rho(B_*)
  &{=}&X^{(n)}_{s*}&{\buildrel {j_2}\over  \rightarrow}&X
\end{matrix}
\end{eqnarray}
where $\pi\colon \mathbb P(j_2^*T_X^*)\rightarrow X^{(n)}_{s*}$
is the natural projection of the $\mathbb P^1$-bundle.
By the Lemma~3.1 in \cite{L-Q}, the restriction of
$R^1(f_{1, 0})_*(ev_1^*T_{\Xn})$ to $\frak U_0$ is
a locally free sheaf of rank $(2d-1)$.
Since the excess dimension of $\frak U_0$ is $(2d-1)$,
Proposition~\ref{virtual-prop} implies that
if $\frak M$ is a closed substack of $\overline{\frak M}_{0, k}
(\Xn, d\beta_n)$ contained in $\frak U_k$, then
\begin{eqnarray}\label{virtual-formula}
[\overline {\frak M}_{0, k}(\Xn, d\beta_n)]^{\text{\rm vir}}|_{\frak M}
= \left \{ {\tilde {f}_{k, 0}}^*\big(c_{2d-1}(R^1(f_{1, 0})_*(ev_{1})^*
T_{\Xn})|_{f_{k, 0}(\frak M)}\big) \right \} |_{\frak M}.
\end{eqnarray}
The following summerizes the formula (32),
Lemma~3.2 and Remark~3.1 in \cite{L-Q}.

\begin{lemma}  \label{lemma_obs}
\begin{enumerate}
\item[{\rm (i)}] $\mathcal O_{B_*}(B_*) \cong
j_1^*\mathcal O_{\mathbb P(j_2^*T_X^*)}(-2)$.

\item[{\rm (ii)}] Let $\mathcal V$ denote the restriction of
$R^1(f_{1, 0})_*(ev_{1})^*T_{\Xn}$ to $\frak U_0$. Then,
the locally free sheaf $\mathcal V$ sits in the exact sequence
\begin{eqnarray*}
0 \to (j_2 \circ \phi)^*\mathcal O_X(-K_X) \to \mathcal V
\to \mathcal E \to 0
\end{eqnarray*}
where $\mathcal E = R^1 (\rf)_*(j_1 \circ \rev)^*((j_2\circ\pi)^*T_{X }
\otimes \mathcal O_{\mathbb P(j_2^*T_X^*)}(-1))$.

\item[{\rm (iii)}] Over $\phi^{-1}(2x_2+x_3+\ldots+x_n)\cong
\overline{\frak M}_{0, 0}(\mathbb P^1, d[\mathbb P^1])$
where $x_2, \ldots, x_n$ are distinct points in $X$,
there is an isomorphism of locally free sheaves:
\begin{eqnarray*}
\mathcal E|_{\phi^{-1}(2x_2+x_3+\ldots+x_n)} \cong
R^1 (f_{1, 0})_*(ev_1)^*(\mathcal O_{\mathbb P^1}(-1) \oplus
\mathcal O_{\mathbb P^1}(-1)).
\end{eqnarray*}
\end{enumerate}
\end{lemma}

Next, using Lemma~\ref{lemma_obs}, we compute other $2$-point
Gromov-Witten invariants.

\begin{lemma} \label{2_obs}
Let $X = \PP^2$ and $d \ge 1$. Then,
\begin{enumerate}
\item[{\rm (i)}]
$\langle \pd(A_1), \pd(A_2) \rangle_{0, d} = 0$ for
the two choices of $(A_1, A_2)$:
\begin{eqnarray*}
(\frak a_{-1}(X)\frak a_{-2}(x) \vac,
  \frak a_{-2}(X)\frak a_{-1}(\ell)\vac),\,\,
(\frak a_{-2}(\ell)\frak a_{-1}(\ell)\vac,
  \frak a_{-1}(X)\frak a_{-2}(\ell)\vac);
\end{eqnarray*}

\item[{\rm (ii)}]
$\langle \pd(A_1), \pd(A_2) \rangle_{0, d} = -4(K_X \cdot \ell)/d$
for the two choices of $(A_1, A_2)$:
\begin{eqnarray*}
(\frak a_{-2}(X)\frak a_{-1}(x) \vac,
  \frak a_{-1}(X)\frak a_{-2}(\ell)\vac),\,\,
(\frak a_{-2}(\ell)\frak a_{-1}(\ell)\vac,
  \frak a_{-2}(X)\frak a_{-1}(\ell)\vac).
\end{eqnarray*}
\end{enumerate}
\end{lemma}
\begin{proof}
(i) Since the proofs for the two choices of $(A_1, A_2)$ are similar,
we only prove $\langle \pd(A_1), \pd(A_2) \rangle_{0, d} = 0$
for $A_1 = \frak a_{-1}(X)\frak a_{-2}(x) \vac$ and
$A_2 = \frak a_{-2}(X)\frak a_{-1}(\ell)\vac$.
Fix a point $x$ and a line $\ell$ in $X = \PP^2$ such that
$x \not \in \ell$. By (\ref{geom-rep}), we see that up to a scalar,
$A_1$ is represented by the closure of the subset
\begin{eqnarray}  \label{2_obs.1}
\{ x'+ \xi|\,\, \xi \in M_2(x) \,\, \text{and } x' \ne x \}.
\end{eqnarray}
Similarly, $A_2$ is represented by the closure of the subset
\begin{eqnarray}  \label{2_obs.2}
\{ \xi + x_1|\,\, x_1 \in \ell, \xi \in M_2(x_2) \,\,
\text{for some } x_2 \not \in \ell\}.
\end{eqnarray}

Working with algebraic cycles instead of cohomology classes, we have
\begin{eqnarray}   \label{2_obs.3}
\langle \pd(A_1), \pd(A_2) \rangle_{0, d} =
[\overline {\frak M}_{0, 2}(X^{[3]}, d\beta_3)]^{\text{\rm vir}}
\cdot ev_2^*[A_1 \times A_2].
\end{eqnarray}
Note that $ev_2^*[A_1 \times A_2]$ is an algebraic cycle supported in
$ev_2^{-1}(A_1 \times A_2)$. By (\ref{2_obs.1}) and (\ref{2_obs.2}),
$ev_2^{-1}(A_1 \times A_2)$ parametrizes all the stable maps
$[\mu: (C; p_1, p_2) \to X^{[3]}]$ satisfying $\rho(\mu(C)) \in 2x + \ell$.
In particular, $ev_2^{-1}(A_1 \times A_2) \subset \frak U_2$.
Applying (\ref{virtual-formula}) to $\frak M = ev_2^{-1}(A_1 \times A_2)$
and combining with Lemma~\ref{lemma_obs}~(ii), we obtain
\begin{eqnarray} \label{2_obs.4}
[\overline {\frak M}_{0, 2}(X^{[3]}, d\beta_3)]^{\text{\rm vir}}|_{\frak M}
&=&\left \{ {\tilde {f}_{2, 0}}^*\big(c_{2d-1}(R^1(f_{1, 0})_*(ev_{1})^*
   T_{\Xn})|_{f_{2, 0}(\frak M)}\big) \right \} |_{\frak M}
   \nonumber  \\
&=&\left \{ {\tilde {f}_{2, 0}}^*\big((j_2 \circ \phi)^*(-K_X) \cdot
   c_{2d-2}(\mathcal E)|_{f_{2, 0}(\frak M)}\big) \right \} |_{\frak M}.
   \qquad
\end{eqnarray}
Now $(j_2 \circ \phi)^*(-K_X)
= 3(j_2 \circ \phi)^*[\ell']$ where the line $\ell'$
in $X = \PP^2$ is chosen not to contain the fixed point $x$.
We have $(j_2 \circ \phi)^{-1}(\ell') \cap
f_{2, 0}(\frak M) = \emptyset$. Therefore,
$(j_2 \circ \phi)^*(-K_X)|_{f_{2, 0}(\frak M)} =0$.
By (\ref{2_obs.4}), $[\overline {\frak M}_{0, 2}
(X^{[3]}, d\beta_3)]^{\text{\rm vir}}|_{\frak M} = 0$.
Since $ev_2^*[A_1 \times A_2]$ is supported in
$\frak M = ev_2^{-1}(A_1 \times A_2)$, we see from
(\ref{2_obs.3}) that $\langle \pd(A_1), \pd(A_2) \rangle_{0, d} = 0$.

(ii) Again, the proofs for the two choices of $(A_1, A_2)$ are similar.
So we only prove $\langle \pd(A_1), \pd(A_2) \rangle_{0, d} =
-4(K_X \cdot \ell)/d$ for $A_1 = \frak a_{-2}(\ell)\frak a_{-1}(\ell)\vac$
and $A_2 = \frak a_{-2}(X)\frak a_{-1}(\ell)\vac$.
We follow the argument for the Lemma~3.3~(ii) in \cite{L-Q}.

Fix three lines $\ell_1, \ell_2, \ell_3 \subset X = \PP^2$ in general position.
Then $A_1$ is represented by the closure of the subset
$\{ \xi +x|\,\, \xi \in M_2(\ell_1), x \in \ell_2, x \not \in |\Supp(\xi)|\}$.
Similarly, $A_2$ is represented by the closure of the subset
\begin{eqnarray*}
\{ \xi+x|\,\, \xi \in M_2(X), x \in \ell_3, x \not \in |\Supp(\xi)|\}.
\end{eqnarray*}
So $ev_2^{-1}(A_1 \times A_2)$ parametrizes all the stable maps
$[\mu: (C; p_1, p_2) \to X^{[3]}]$ satisfying
$\rho(\mu(C)) \in 2\ell_1 + (\ell_2 \cap \ell_3) \subset B_*$,
and $ev_2^*[A_1 \times A_2]$ is a cycle in
$ev_2^{-1}(A_1 \times A_2) \subset \frak U_2$.
As in (\ref{2_obs.3}) and (\ref{2_obs.4}), we see that
$\langle \pd(A_1), \pd(A_2) \rangle_{0, d}$ is equal to
\begin{eqnarray*}
& &{\tilde {f}_{2, 0}}^*\big((j_2 \circ \phi)^*(-K_X) \cdot
   c_{2d-2}(\mathcal E)\big) \,\, ev_2^*[A_1 \times A_2].
\end{eqnarray*}
Since ${\tilde {f}_{2, 0}}^*\big((j_2 \circ \phi)^*(-K_X) \cdot
c_{2d-2}(\mathcal E)\big)$ is supported in $\frak U_2$,
recalling the definition of $\tilde {ev}_2$ from the paragraph
containing (\ref{com-diagram2}), we see that
$\langle \pd(A_1), \pd(A_2) \rangle_{0, d}$ equals
\begin{eqnarray*}
{\tilde {f}_{2, 0}}^*\big((j_2 \circ \phi)^*(-K_X) \cdot
   c_{2d-2}(\mathcal E)\big) \,\,
   {\tilde {ev}_2}^*\big ( ([A_1] [B_*]) \times ([A_2] [B_*]) \big ).
\end{eqnarray*}
Now, $[A_i] [B_*] = [A_i \cap B_*] c_1(\mathcal O_{B_*}(B_*))$.
Let $\mathbb D$ stand for the first Chern class of the tautological
line bundle over $B_* \cong \mathbb P(j_2^*T_X^*)$.
Then we obtain from Lemma~\ref{lemma_obs}~(i) that the invariant
$\langle \pd(A_1), \pd(A_2) \rangle_{0, d}$ is equal to
\begin{eqnarray}   \label{2_obs.5}
4{\tilde {f}_{2, 0}}^*\big((j_2 \circ \phi)^*(-K_X) \cdot
   c_{2d-2}(\mathcal E)\big) \cdot
   {\tilde {ev}_2}^*\big ( ([A_1 \cap B_*] \mathbb D) \times
   ([A_2 \cap B_*] \mathbb D) \big ).
\end{eqnarray}

Fix a line $\ell$ such that $\ell_1, \ell_2, \ell_3, \ell$ are
in general position. We claim that
\begin{eqnarray}   \label{2_obs.6}
{\tilde {f}_{2, 0}}^*(j_2 \circ \phi)^*[\ell] \cdot
   {\tilde {ev}_2}^*\big ( ([A_1 \cap B_*] \mathbb D) \times
   ([A_2 \cap B_*] \mathbb D) \big )
= [{\tilde {ev}_2}^{-1}(\xi_1 \times \xi_2)]
\end{eqnarray}
where $\xi_1$ and $\xi_2$ are two fixed points in $M_2(x_1)+x_2$ with
$\{x_1 \}= \ell_1 \cap \ell$, and $\{x_2 \}= \ell_2 \cap \ell_3$.
To see this, let ${\tilde e}_1$ and ${\tilde e}_2$ be the restrictions
to $\frak U_2$ of the two evaluation maps from
$\overline {\frak M}_{0, 2}(X^{[3]}, d\beta_3)$ to $X^{[3]}$.
We regard ${\tilde e}_1$ and ${\tilde e}_2$ as morphisms
from $\frak U_2$ to $B_*$.
Then, ${\tilde {ev}_2} = {\tilde e}_1 \times {\tilde e}_2$ and
$\phi \circ {\tilde {f}_{2, 0}} = \rho \circ {\tilde {e}_1}$. So
\begin{eqnarray*}
& &{\tilde {f}_{2, 0}}^*(j_2 \circ \phi)^*[\ell] \cdot
   {\tilde {ev}_2}^*\big ( ([A_1 \cap B_*] \mathbb D) \times
   ([A_2 \cap B_*] \mathbb D) \big )   \\
&=&{\tilde {f}_{2, 0}}^*(j_2 \circ \phi)^*[\ell] \cdot
   {\tilde {e}_1}^*([A_1 \cap B_*] \mathbb D)  \cdot
   {\tilde {e}_2}^*([A_2 \cap B_*] \mathbb D) \\
&=&{\tilde {e}_1}^*((j_2 \circ \rho)^*[\ell] \cdot [A_1 \cap B_*]
   \cdot \mathbb D) \cdot {\tilde {e}_2}^*([A_2 \cap B_*] \mathbb D).
\end{eqnarray*}
Now the cycle $(j_2 \circ \rho)^*[\ell] \cdot [A_1 \cap B_*] \cdot
\mathbb D$ is represented by $\eta_1 + \ell_2$ where
$\eta_1$ is a fixed point in $M_2(x_1)$.
So ${\tilde {e}_1}^*((j_2 \circ \rho)^*[\ell] \cdot [A_1 \cap B_*]
\cdot \mathbb D)$ is represented by the substack $\frak M_2$ of
$\overline {\frak M}_{0, 2}(X^{[3]}, d\beta_3)$
parametrizing all the stable maps $[\mu: (C; p_1, p_2) \to X^{[3]}]$
such that $\mu(C) = M_2(x_1) + x$ for some $x \in \ell_2$ and
$\mu(p_1) = \eta_1 + x$. It follows that
\begin{eqnarray*}
& &{\tilde {f}_{2, 0}}^*(j_2 \circ \phi)^*[\ell] \cdot
   {\tilde {ev}_2}^*\big ( ([A_1 \cap B_*] \mathbb D) \times
   ([A_2 \cap B_*] \mathbb D) \big )   \\
&=&[\frak M_2] \cdot {\tilde {e}_2}^*([A_2 \cap B_*] \mathbb D)
   = [{\tilde {ev}_2}^{-1}(\xi_1 \times \xi_2)]
\end{eqnarray*}
where $\xi_1 = \eta_1 + x_2$ and $\xi_2$ is a fixed point in $M_2(x_1)+x_2$.
This proves (\ref{2_obs.6}).

By (\ref{2_obs.5}) and (\ref{2_obs.6}), $\langle \pd(A_1), \pd(A_2)
\rangle_{0, d}$ is equal to
\begin{eqnarray}   \label{2_obs.7}
& &12{\tilde {f}_{2, 0}}^*\big (c_{2d-2}(\mathcal E)\big) \cdot
   [{\tilde {ev}_2}^{-1}(\xi_1 \times \xi_2)]
   \nonumber  \\
&=&-4(K_X \cdot \ell) \cdot c_{2d-2}(\mathcal E) \cdot ({\tilde {f}_{2, 0}})_*
   [{\tilde {ev}_2}^{-1}(\xi_1 \times \xi_2)].
\end{eqnarray}
Note that ${\tilde {ev}_2}^{-1}(\xi_1 \times \xi_2)$ parametrizes
all the stable maps $[\mu: (C; p_1, p_2) \to X^{[3]}]$
in $\overline {\frak M}_{0, 2}(X^{[3]}, d\beta_3)$
satisfying $\mu(p_1) = \xi_1$ and $\mu(p_2) = \xi_2$.
For these stable maps, we must have $\mu(C) = M_2(x_1)+x_2$.
So the restriction of ${\tilde {f}_{2, 0}}$ to
${\tilde {ev}_2}^{-1}(\xi_1 \times \xi_2)$ is a degree-$d^2$ morphism
to $\phi^{-1}(2x_1+x_2)$. Thus, $({\tilde {f}_{2, 0}})_*
[{\tilde {ev}_2}^{-1}(\xi_1 \times \xi_2)] = d^2[\phi^{-1}(2x_1+x_2)]$.
By (\ref{2_obs.7}), we obtain $\langle \pd(A_1), \pd(A_2) \rangle_{0, d}
= -4(K_X \cdot \ell)d^2 \cdot c_{2d-2}(\mathcal E|_{\phi^{-1}(2x_1+x_2)})$.
By Lemma~\ref{lemma_obs}~(iii) and the Theorem 9.2.3 in \cite{C-K},
$c_{2d-2}(\mathcal E|_{\phi^{-1}(2x_1+x_2)}) = 1/d^3$. Therefore, we have
$\langle \pd(A_1), \pd(A_2) \rangle_{0, d} = -4(K_X \cdot \ell)/d$.
\end{proof}

In view of Lemma~\ref{2_zero}, Lemma~\ref{2_zero_K3} and Lemma~\ref{2_obs},
the only $2$-point Gromov-Witten invariant
$\langle \pd(A_1), \pd(A_2) \rangle_{0, d}$
with $A_1 \in \frak B_6$ and $A_2 \in \frak B_8$ that has not been computed
is when $A_1 = \frak a_{-3}(\ell)\vac$ and $A_2 = \frak a_{-3}(X)\vac$.
This invariant
\begin{eqnarray}  \label{2_remain}
\langle \pd(\frak a_{-3}(\ell)\vac),
\pd(\frak a_{-3}(X)\vac) \rangle_{0, d}
\end{eqnarray}
will be studied in Sect.~\ref{sect_loc} by using the localization formula.

We summarize the results in this subsection into a theorem.

\begin{theorem}  \label{thm_2_point}
Let $X = \PP^2$, and $\frak B_6$ and $\frak B_8$ be defined in
Definition~\ref{def_B_i}.
Let $d \ge 1$, $A_1 \in \frak B_6$ and $A_2 \in \frak B_8$.
Let $x, \ell$ be a point and a line in $X$ respectively.
Then, $\langle \pd(A_1), \pd(A_2) \rangle_{0, d}$ is zero unless
the pair $(A_1, A_2)$ is one of the following:
\par
{\rm (i)}
$(\mathfrak a_{-2}(X)\mathfrak a_{-1}(x) \vac,
\mathfrak a_{-1}(X)\mathfrak a_{-2}(\ell)\vac)$
\par
{\rm (ii)} $(\mathfrak a_{-2}(\ell)\mathfrak a_{-1}(\ell)\vac,
\mathfrak a_{-2}(X)\mathfrak a_{-1}(\ell)\vac)$
\par
{\rm (iii)} $(\mathfrak a_{-3}(\ell)\vac, \mathfrak a_{-3}(X)\vac)$.
\par\noindent
Moreover, $\langle \pd(A_1), \pd(A_2) \rangle_{0, d} =
12/d$ in cases (i) and (ii). \qed
\end{theorem}

\subsection{$\langle \pd(A_1), \pd(A_2), \pd(A_3) \rangle_{0, d}$
with $A_1, A_2, A_3 \in \frak B_8$}
\label{subsect_3}
\par
$\,$

\begin{lemma} \label{3_zero}
The Gromov-Witten invariants
$\langle \pd(A_1), \pd(A_2), \pd(A_3) \rangle_{0, d}$ are equal to
zero for the following triples of $(A_1, A_2, A_3) \in (\frak B_8)^3$:
\begin{eqnarray*}
&A_1= \frak a_{-3}(X)\vac, \,\,A_2 \ne \frak a_{-3}(X)\vac, \,\,
  A_3 \ne \frak a_{-3}(X)\vac, & \\
&A_1= A_2 = \frak a_{-1}(X)\frak a_{-1}(\ell)^2\vac, \,\,
  A_3 \,\, \text{\rm arbitrary},& \\
&A_1= \frak a_{-1}(X)^2\frak a_{-1}(x)\vac, \,\,
  A_2 \,\, \text{\rm arbitrary}, \,\, A_3 \,\, \text{\rm arbitrary}, & \\
&(\frak a_{-3}(X)\vac, \frak a_{-3}(X)\vac,
  \frak a_{-1}(X)\frak a_{-1}(\ell)^2\vac), & \\
&A_1=A_2=\frak a_{-2}(X)\frak a_{-1}(\ell)\vac, \,\,
  A_3 \ne \frak a_{-1}(X)\frak a_{-2}(\ell)\vac, & \\
&(\frak a_{-2}(X)\frak a_{-1}(\ell)\vac,
  \frak a_{-1}(X)\frak a_{-2}(\ell)\vac,
  \frak a_{-1}(X)\frak a_{-1}(\ell)^2\vac), & \\
&A_1, A_2, A_3 \in \{\frak a_{-1}(X)\frak a_{-2}(\ell)\vac,
  \frak a_{-1}(X)\frak a_{-1}(\ell)^2\vac\}.&
\end{eqnarray*}
\end{lemma}
\begin{proof}
The arguments are similar to those for Lemma~\ref{2_zero} and
Lemma~\ref{2_zero_K3}.
\end{proof}

\begin{lemma} \label{3_obs}
Let $X = \PP^2$, $\ell \subset X$ be a line, and $d \ge 1$. Then,
\begin{enumerate}
\item[{\rm (i)}]
$\langle \pd(A_1), \pd(A_2), \pd(A_3) \rangle_{0, d} = 0$
for the following triple:
\begin{eqnarray*}
(A_1, A_2, A_3) = (\frak a_{-2}(X)\frak a_{-1}(\ell)\vac,
\frak a_{-1}(X)\frak a_{-2}(\ell)\vac,
\frak a_{-1}(X)\frak a_{-2}(\ell)\vac);
\end{eqnarray*}

\item[{\rm (ii)}]
$\langle \pd(A_1), \pd(A_2), \pd(A_3) \rangle_{0, d} = 8(K_X \cdot \ell)$
for the triple:
\begin{eqnarray*}
(A_1, A_2, A_3)= (\frak a_{-2}(X)\frak a_{-1}(\ell)\vac,
\frak a_{-2}(X)\frak a_{-1}(\ell)\vac,
\frak a_{-1}(X)\frak a_{-2}(\ell)\vac).
\end{eqnarray*}
\end{enumerate}
\end{lemma}
\begin{proof}
The arguments are similar to those for Lemma~\ref{2_obs}~(i) and (ii).
\end{proof}

According to Lemma~\ref{3_zero} and Lemma~\ref{3_obs},
it remains to compute the invariants
$\langle \pd(A_1), \pd(A_2), \pd(A_3) \rangle_{0, d}$
for the following 3 triples of $(A_1, A_2, A_3) \in (\frak B_8)^3$:
\begin{eqnarray*}
&A_1= A_2 = \frak a_{-3}(X)\vac,& \\
&A_3 = \frak a_{-2}(X)\frak a_{-1}(\ell)\vac,
  \frak a_{-1}(X)\frak a_{-2}(\ell)\vac, \frak a_{-3}(X)\vac.&
\end{eqnarray*}
In the next two lemmas, we shall calculate them in terms
of (\ref{2_remain}). Put
\begin{eqnarray}  \label{E_i}
\mathcal E_i = \pi_1(\pi_2^*\mathcal O_X(i)|_{\mathcal O_{\mathcal Z_3}})
\end{eqnarray}
where $\pi_1$ and $\pi_2$
denote the projections of $X^{[3]} \times X$ to the two factors.
It is known that $c_1(\mathcal E_i) = i D_\ell - B_3/2$.
Using the commutation relations among standard operators on $\mathbb H$
(e.g. the Theorem 3.1 in \cite{LQW4}), we obtain
\begin{eqnarray}   \label{a_3_X}
   c_1(\mathcal E_0)^2
&=&\frak a_{-3}(X)\vac
   -\frak a_{-1}(X)^2\frak a_{-1}(x)\vac \nonumber \\
& &- {1 \over 2}\frak a_{-1}(X)\frak a_{-1}(\ell)^2\vac
   - {1 \over 2}\frak a_{-1}(X)\frak a_{-2}(K_X)\vac.
\end{eqnarray}
%\begin{eqnarray}   \label{a_3_X}
%   c_1(\mathcal E_0)^2
%&=&\frak a_{-3}(X)\vac
%   -\frak a_{-1}(X)^2\frak a_{-1}(x)\vac \nonumber \\
%& &- {1 \over 2}\frak a_{-1}(X)\frak a_{-1}(\ell)^2\vac
%   + {3 \over 2}\frak a_{-1}(X)\frak a_{-2}(\ell)\vac.
%\end{eqnarray}

\begin{lemma} \label{3_comp_law1}
Let $d \ge 1$ and $A = \frak a_{-3}(X)\vac$.
Let $w_1, w_2$ denote the two invariants
$\langle \pd(A), \pd(A), \pd(A_3) \rangle_{0, d}$
for $A_3 = \frak a_{-2}(X)\frak a_{-1}(\ell)\vac, \,\,
\frak a_{-1}(X)\frak a_{-2}(\ell)\vac$ respectively.
Then, $w_1 = w_2 = -2d \,\, \langle \pd(\frak a_{-3}(\ell)\vac),
\pd(\frak a_{-3}(X)\vac) \rangle_{0, d}$.
\end{lemma}
\begin{proof}
Since the arguments for $w_1$ and $w_2$ are almost the same,
we only prove that $w_2 = -2d \,\, \langle
\pd(\frak a_{-3}(\ell)\vac), \pd(\frak a_{-3}(X)\vac) \rangle_{0, d}$.
Let $c_1 = c_1(\mathcal E_0) = -B_3/2$ (we regard a divisor as
either a homology class or a cohomology class depending on the context).
Apply the composition law (\ref{comp_law}) to
$\alpha_1 = \alpha_2 = c_1, \alpha_3 =\pd(\frak a_{-3}(X)\vac),
\alpha_4 = \pd(\frak a_{-1}(X)\frak a_{-2}(\ell)\vac)$,
and to the basis $\{ \Delta_a \}$ of $H^*(X^{[3]})$ given by (\ref{basis}).

First of all, the left-hand-side of (\ref{comp_law}) is equal to
\begin{eqnarray} \label{3_comp_law1.1}
& &\langle c_1^2, \alpha_3, \alpha_4 \rangle_{0, d}
   +\langle c_1, c_1, \alpha_3 \alpha_4 \rangle_{0, d}
   \nonumber \\
&+&\sum_{d_1 + d_2 = d, \, d_1, d_2 > 0} \,\,
   \sum_{a} \,\, \langle c_1, c_1, \Delta_a \rangle_{0, d_1}
   \, \langle \Delta^a, \alpha_3, \alpha_4 \rangle_{0, d_2}.
\end{eqnarray}
By (\ref{a_3_X}) and Lemma~\ref{3_zero},
$\langle c_1^2, \alpha_3, \alpha_4 \rangle_{0, d} = w_2$.
Since the intersection number $(c_1 \cdot \beta_3)$ is equal to $1$,
$\langle c_1, c_1, \alpha_3 \alpha_4 \rangle_{0, d} =
d^2 \, \langle \alpha_3 \alpha_4 \rangle_{0, d}$ and
$\langle c_1, c_1, \Delta_a \rangle_{0, d_1} =
d_1^2 \, \langle \Delta_a \rangle_{0, d_1}$.
By Lemma~\ref{1pt}, $\langle \Delta_a \rangle_{0, d_1} \ne 0$
only when $\Delta_a = \pd(\frak a_{-2}(\ell)\frak a_{-1}(x)\vac)$.
Note that $\Delta^a = -1/2 \, \pd(\frak a_{-1}(X)\frak a_{-2}(\ell)\vac)$.
So $\langle \Delta^a, \alpha_3, \alpha_4 \rangle_{0, d_2} = 0$ by
Lemma~\ref{3_zero}. It follows from (\ref{3_comp_law1.1}) that
the left-hand-side of (\ref{comp_law}) is equal to
\begin{eqnarray} \label{3_comp_law1.2}
w_2 + d^2 \, \langle \alpha_3 \alpha_4 \rangle_{0, d}.
\end{eqnarray}

We claim that $\langle \alpha_3 \alpha_4 \rangle_{0, d}
= -12(K_X \cdot \ell)/d^2$.
To prove this, note from (\ref{a_3_X}) that $\frak a_{-3}(X)\vac =
c_1^2 + \frak a_{-1}(X)^2\frak a_{-1}(x)\vac
+ {1/2} \, \frak a_{-1}(X)\frak a_{-1}(\ell)^2\vac
- {3/2} \, \frak a_{-1}(X)\frak a_{-2}(\ell)\vac$.
Choose lines $\ell', \ell''$ in $X = \PP^2$ such that
$\ell, \ell', \ell''$ are in general position.
Then, $\big ( \frak a_{-1}(X)\frak a_{-2}(\ell)\vac \big )
\cap \big ( \frak a_{-1}(X)\frak a_{-2}(\ell')\vac \big )
\cap \big (\frak a_{-1}(X)\frak a_{-2}(\ell'')\vac \big ) = \emptyset$.
It follows that $\big ( \frak a_{-1}(X)\frak a_{-2}(\ell)\vac \big )^3 = 0$.
In view of the linear basis in Lemma~\ref{lemmacohoclass}~(ii),
we see that $\big ( \frak a_{-1}(X)\frak a_{-2}(\ell)\vac \big )^2$
is a linear combination of $\frak a_{-1}(X)\frak a_{-1}(x)^{2} \vac$,
$\frak a_{-1}(\ell)^2 \frak a_{-1}(x)\vac$,
$\frak a_{-1}(\ell)\frak a_{-2}(x)\vac$, and $\frak a_{-3}(x)\vac$.
Hence $\langle \pd(\frak a_{-1}(X)\frak a_{-2}(\ell)\vac)
\,\, \alpha_4 \rangle_{0, d} = 0$ according to Lemma~\ref{1pt}~(i),
and we see that $\langle \alpha_3 \alpha_4 \rangle_{0, d}$ is equal to
\begin{eqnarray*}
\langle c_1^2 \alpha_4 \rangle_{0, d}
+ \langle \pd(\frak a_{-1}(X)^2\frak a_{-1}(x)\vac)
\,\, \alpha_4 \rangle_{0, d}
+ {1 \over 2} \langle \pd(\frak a_{-1}(X)\frak a_{-1}(\ell)^2\vac)
\,\, \alpha_4 \rangle_{0, d}.
\end{eqnarray*}
Since $(D_\ell)^2 = \frak a_{-1}(X)\frak a_{-1}(\ell)^2\vac
+1/2 \, \frak a_{-1}(X)^2\frak a_{-1}(x)\vac$, we obtain
\begin{eqnarray} \label{3_comp_law1.3}
\langle \alpha_3 \alpha_4 \rangle_{0, d}
= \langle c_1^2 \alpha_4 \rangle_{0, d}
+ {1 \over 2} \langle D_\ell^2 \alpha_4 \rangle_{0, d}
+ {3 \over 4} \langle \pd(\frak a_{-1}(X)^2\frak a_{-1}(x)\vac)
\,\, \alpha_4 \rangle_{0, d}.
\end{eqnarray}
Since $\frak a_{-1}(X)^2\frak a_{-1}(x)\vac
\cdot \frak a_{-1}(X)\frak a_{-2}(\ell)\vac
= 2\frak a_{-2}(\ell)\frak a_{-1}(x)\vac$, the third term in
(\ref{3_comp_law1.3}) is equal to $3(K_X \cdot \ell)/d^2$
by Lemma~\ref{1pt}~(ii).
Since $D_\ell^2 \cdot \frak a_{-1}(X)\frak a_{-2}(\ell)\vac
= \frak a_{-2}(\ell)\frak a_{-1}(x)\vac
+ 4\frak a_{-1}(\ell)\frak a_{-2}(x)\vac$, the second term in
(\ref{3_comp_law1.3}) is equal to $(K_X \cdot \ell)/d^2$
by Lemma~\ref{1pt}. Using a similar argument, we see that
the first term in (\ref{3_comp_law1.3}) is equal to
$-16(K_X \cdot \ell)/d^2$.
Thus, $\langle \alpha_3 \alpha_4 \rangle_{0, d} =
-12(K_X \cdot \ell)/d^2$ in view of (\ref{3_comp_law1.3}).

Combining with (\ref{3_comp_law1.2}), we see that the left-hand-side
of (\ref{comp_law}) is equal to $w_2 - 12(K_X \cdot \ell)$.
Similarly, the right-hand-side of (\ref{comp_law}) is equal to
\begin{eqnarray*}
-2d \,\, \langle \pd(\frak a_{-3}(\ell)\vac),
\pd(\frak a_{-3}(X)\vac) \rangle_{0, d} -12 (K_X \cdot \ell).
\end{eqnarray*}
Hence we have $w_2 = -2d \,\, \langle \pd(\frak a_{-3}(\ell)\vac),
\pd(\frak a_{-3}(X)\vac) \rangle_{0, d}$.
\end{proof}

\begin{lemma} \label{3_comp_law2}
Let $d \ge 1$. Put $f(d) = d \, \langle \pd(\frak a_{-3}(\ell)\vac),
\pd(\frak a_{-3}(X)\vac) \rangle_{0, d}$. Let $w_3$ denote
$\langle \pd(A), \pd(A), \pd(A) \rangle_{0, d}$
for $A = \frak a_{-3}(X)\vac$. Then $w_3$ equals
\begin{eqnarray*}
&\displaystyle{-24K_X^2 - 18(K_X \cdot \ell)+5(K_X \cdot \ell)f(d)}& \\
&\displaystyle{-2(K_X \cdot \ell) \sum\limits_{0 < d_1 < d} f(d_1)
+{1 \over 3} \sum\limits_{0 < d_1 < d} f(d_1)f(d-d_1).}&
\end{eqnarray*}
\end{lemma}
\begin{proof}
Our idea is the same as in the proof of Lemma~\ref{3_comp_law1}.
Let $c_1 = c_1(\mathcal E_0)$. Apply (\ref{comp_law})
to $\alpha_1 = \alpha_2 = c_1$ and
$\alpha_3 = \alpha_4 = \pd(\frak a_{-3}(X)\vac)$.
Then, the left-hand-side of (\ref{comp_law}) is still of the form
(\ref{3_comp_law1.1}).
By (\ref{a_3_X}), Lemma~\ref{3_zero} and Lemma~\ref{3_comp_law1},
$\langle c_1^2, \alpha_3, \alpha_4 \rangle_{0, d} =
w_3 - (K_X \cdot \ell)/2 \, w_2 = w_3 + (K_X \cdot \ell)f(d)$.
Also, $\langle c_1, c_1, \alpha_3 \alpha_4 \rangle_{0, d}
= d^2 \, \langle \alpha_3 \alpha_4 \rangle_{0, d}
= 24K_X^2 + 18(K_X \cdot \ell)$, and
$\sum_{a} \,\, \langle c_1, c_1, \Delta_a \rangle_{0, d_1}
   \, \langle \Delta^a, \alpha_3, \alpha_4 \rangle_{0, d_2}$ is equal to
\begin{eqnarray*}
& &-{d_1^2 \over 2} \,\, \langle
   \pd(\frak a_{-2}(\ell)\frak a_{-1}(x)\vac) \rangle_{0, d_1}
   \, \langle \pd(\frak a_{-1}(X) \frak a_{-2}(\ell)\vac),
   \alpha_3, \alpha_4 \rangle_{0, d_2} \\
&=&-{d_1^2 \over 2} \cdot {2(K_X \cdot \ell) \over d_1^2} \cdot
   (-2f(d_2)) = 2(K_X \cdot \ell)f(d_2)
\end{eqnarray*}
by Lemma~\ref{1pt}~(ii) and Lemma~\ref{3_comp_law1}.
So the left-hand-side of (\ref{comp_law}) is
\begin{eqnarray} \label{3_comp_law2.1}
w_3 + (K_X \cdot \ell)f(d) + 24K_X^2 + 18(K_X \cdot \ell) +
2(K_X \cdot \ell) \sum_{0 < d_1 < d} f(d_1).
\end{eqnarray}

Similarly, the right-hand-side of (\ref{comp_law}) is equal to
\begin{eqnarray} \label{3_comp_law2.2}
6(K_X \cdot \ell) f(d) +{1 \over 3} \sum_{0 < d_1 < d} f(d_1)f(d-d_1).
\end{eqnarray}
Now we prove the lemma by comparing (\ref{3_comp_law2.1}) and
(\ref{3_comp_law2.2}).
\end{proof}

The results in this subsection are summarized into a theorem.

\begin{theorem}  \label{thm_3_point}
Let $X = \PP^2$, and $\frak B_8$ be defined in Definition~\ref{def_B_i}.
Let $\ell \subset X$ be a line.
Let $d \ge 1$, $f(d) = d \, \langle \pd(\frak a_{-3}(\ell)\vac),
\pd(\frak a_{-3}(X)\vac) \rangle_{0, d}$,
and $A_1, A_2, A_3 \in \frak B_8$.
Then, the $3$-point genus-$0$ Gromov-Witten invariant
$\langle \pd(A_1), \pd(A_2), \pd(A_3) \rangle_{0, d}$ is zero unless
the unordered triple $(A_1, A_2, A_3)$ is one of the following:
\par
{\rm (i)}
$(\frak a_{-2}(X)\frak a_{-1}(\ell)\vac,
  \frak a_{-2}(X)\frak a_{-1}(\ell)\vac,
  \frak a_{-1}(X)\frak a_{-2}(\ell)\vac)$
\par
{\rm (ii)} $(\mathfrak a_{-3}(X)\vac, \mathfrak a_{-3}(X)\vac,
\frak a_{-2}(X)\frak a_{-1}(\ell)\vac)$
\par
{\rm (iii)} $(\mathfrak a_{-3}(X)\vac, \mathfrak a_{-3}(X)\vac,
\frak a_{-1}(X)\frak a_{-2}(\ell)\vac)$
\par
{\rm (iv)} $(\mathfrak a_{-3}(X)\vac, \mathfrak a_{-3}(X)\vac,
\mathfrak a_{-3}(X)\vac)$.
\par\noindent
Moreover, $\langle \pd(A_1), \pd(A_2), \pd(A_3) \rangle_{0, d}
= -24$ for case (i); for cases (ii) and (iii),
$\langle \pd(A_1), \pd(A_2), \pd(A_3) \rangle_{0, d} = -2f(d)$;
for case (iv),
\begin{eqnarray*}
& &\langle \pd(A_1), \pd(A_2), \pd(A_3) \rangle_{0, d} \\
&=&-162-15f(d) +6\sum\limits_{0 < d_1 < d} f(d_1)
+{1 \over 3} \sum\limits_{0 < d_1 < d} f(d_1)f(d-d_1). \qed
\end{eqnarray*}
\end{theorem}

\section{\bf Computation of $\langle \pd(\frak a_{-3}(\ell)\vac),
\pd(\frak a_{-3}(X)\vac) \rangle_{0, d}$}
\label{sect_loc}

In this section, we study the remaining $2$-point Gromov-Witten invariant
\begin{eqnarray*}
\langle \pd(\frak a_{-3}(\ell)\vac),
\pd(\frak a_{-3}(X)\vac) \rangle_{0, d}
\end{eqnarray*}
in (\ref{2_remain}). Using the standard $(\C^*)^2$-action on $X = \PP^2$
and the virtual localization formula in \cite{G-P},
we reduce the computation to a summation over stable graphs.
This allows us to calculate $\langle \pd(\frak a_{-3}(\ell)\vac),
\pd(\frak a_{-3}(X)\vac) \rangle_{0, d}$ for $d \le 4$.

\subsection{The contracted $(\C^*)^2$-invariant curves in $(\PP^2)^{[3]}$}
\label{subsect_curve}
$\,$

Let $T \subset \SL_3(\C)$ be the subgroup consisting of diagonal matrices.
Then $T \simeq (\C^*)^2$ acts on $\PP^2$ with fixed
points $P_0 = (1, 0, 0)$, $P_1 = (0, 1, 0)$ and $P_2 = (0, 0, 1)$.
There is an induced action of $T$ on the Hilbert scheme
$(\PP^2)^{[3]}$ with a finite number of fixed points.
The $T$-fixed points in $(\PP^2)^{[3]}$ are enumerated as follows.
If $(u_i, v_i)$ are the local coordinates at the fixed point $P_i$,
then there are three $T$-fixed points
%$Q_{i, 1}$, $Q_{i, 2}$, $Q_{i, 3}$
in $M_3({P_i}) \subset (\PP^2)^{[3]}$ corresponding to
the partitions $(3)$, $(2,1)$ and $(1,1,1)$ of $3$.
The corresponding ideals are $(u_i^3, v_i)$, $(u_i^2, u_iv_i, v_i^2)$
and $(u_i,v_i^3)$. Also for each ordered pair of points $(P_i,P_j)$
with $i \ne j$, we have two fixed points $R_{i,j}^{(1)} = \xi_{i, 1} + P_j$
and $R_{i,j}^{(2)} = \xi_{i, 2} + P_j$ in $(\PP^2)^{[3]}$,
where $\xi_{i, 1}, \xi_{i, 2} \in M_2(P_i)$ correspond to the ideals
$(u_i, v_i^2), (u_i^2, v_i)$ respectively.
Finally, $P_0+P_1+P_2$ is also a $T$-fixed point in $(\PP^2)^{[3]}$.

Next, we start enumerating $T$-invariant curves. Observe that
a $T$-invariant curve is the closure of a 1-dimensional $T$-orbit.
Thus, a $T$-invariant curve is the $T$-orbit of a point in
a fixed component of a 1-parameter subgroup of $T$ corresponding
to the kernel of the $T$-action along the curve.
In particular a $T$-invariant curve is a smooth rational curve,
and must contain exactly two fixed points.

We are only interested in $T$-invariant curves that are contracted
under the Hilbert-Chow morphism $(\PP^2)^{[3]} \to (\PP^2)^{(3)}$.
Such curves must be entirely contained in $M_3({P_i})$ for some $i$,
or in $M_2({P_i}) + P_j$ for some $i \ne j$.
Since $M_2({P_i}) \simeq \PP^1$, we immediately obtained
six $T$-invariant curves $C_{i,j} \, {\overset {\rm def} =} \,
M_2({P_i}) + P_j$, with $1 \le i, j \le 3$ and $i \ne j$, contracted
by the Hilbert-Chow morphism $(\PP^2)^{[3]} \to (\PP^2)^{(3)}$.

We now analyze $T$-invariant curves in
$M_3(P_i)$, by using a tangent space analysis.
%First of all, suppose that $Z \subset (\PP^2)^{[3]}$ is
%a component of the fixed locus of a $1$-parameter subgroup $G$.
%A more or less explicit argument shows that if $Z$ contains any curve
%contracted to $3P_i$, then $Z$ is entirely contained in $M_3(P_i)$.
%Next, we use a tangent space analysis to determine the
%$T$-invariant curves in $M_3(P_i)$.
Suppose that $(s,t)(u_i, v_i) = (\lambda_i(s,t)u_i, \mu_i(s,t)v_i)$
where $\lambda_i$ and $\mu_i$ are independent characters of $T$.
Let $Q_{i, 0}, Q_{i, 1}, Q_{i, 2} \in M_3(P_i)$
be the three $T$-fixed points corresponding to the ideals
$(u_i^2,u_iv_i,v_i^2), (u_i^3,v_i), (u_i,v_i^3)$ respectively.
For simplicity, denote the tangent space of $(\PP^2)^{[3]}$
at the point $Q_{i, j}$ by $T_{Q_{i, j}}$.
By \cite{E-S}, we have the following decompositions
for the tangent spaces as a representation of $T$:
\begin{eqnarray}
T_{Q_{i, 0}}
&=&2\lambda_i^{-1} + 2\mu_i^{-1}
   + \lambda_i^{-2}\mu_i + \lambda_i \mu_i^{-2}   \label{eq.tq0}  \\
T_{Q_{i, 1}}
&=&\lambda_i^{-1}\mu_i^2 + \lambda_i^{-1}\mu_i+ \lambda_i^{-1} +
   \mu_i^{-3} + \mu_i^{-2} + \mu_i^{-1}  \label{eq.tq1}  \\
T_{Q_{i, 2}}
&=&\lambda_i^{-3} + \lambda_i^{-2} + \lambda_i^{-1} + \lambda_i^2\mu_i^{-1}
  + \lambda_i \mu_i^{-1} + \mu_i^{-1}. \label{eq.tq2}
\end{eqnarray}
The kernel of each character appearing in equations (\ref{eq.tq0}),
(\ref{eq.tq1}), (\ref{eq.tq2}) determines $1$-parameter subgroup
whose fixed locus contains $T$-invariant curves.
Since we are interested only in $T$-invariant curves
contained in $M_3(P_i)$, we need only to analyze characters of
the form $\lambda_i^k \mu_i^\ell$ with $k \ell \neq 0$.
(The kernel of a character $\lambda_i^k$ or $\mu_i^\ell$ will
have fixed locus that moves out of the punctual Hilbert scheme.)

Looking at $T_{Q_{i, 0}}$ we see that the character
$\lambda_i \mu_i^{-2}$ has multiplicity one.
This means that its kernel has one-dimensional
fixed component containing the point $Q_{i, 0}$. Now the character
$\lambda_i^{-1}\mu_i^2$ in $T_{Q_{i, 1}}$ has the same kernel as
the character $\lambda_i \mu_i^{-2}$ in $T_{Q_{i, 0}}$.
So there is a unique $T$-invariant curve, denoted by $C_{0, 1}^{(i)}$,
which contains $Q_{i, 0}$ and $Q_{i, 1}$, and is the fixed locus of
$\ker(\lambda_i \mu_i^{-2})$.
Similar analysis shows that there are two other $T$-invariant curves
$C_{0, 2}^{(i)}$ and $C_{1, 2}^{(i)}$ in $M_3(P_i)$;
namely, $C_{0, 2}^{(i)}$ through $Q_{i, 0}$ and $Q_{i, 2}$ which is
the fixed locus of $\ker(\lambda_i^{-2}\mu_i)$,
while $C_{1, 2}^{(i)}$ through $Q_{i, 1}$ and $Q_{i, 2}$ which is
the fixed locus of $\ker(\lambda_i^{-1}\mu_i)$.
This analysis partially proves the following.

\begin{lemma} \label{15_curve}
There are 15 $T$-invariant curves contracted under the Hilbert-Chow
morphism $(\PP^2)^{[3]} \to (\PP^2)^{(3)}$. They are described as follows:
\begin{enumerate}
\item[{\rm (i)}] the six curves $C_{i, j} = M_2(P_i) + P_j$
where $1 \le i, j \le 3$ and $i \neq j$;

\item[{\rm (ii)}] the nine curves $C_{k, \ell}^{(i)} \subset M_3(P_i)$
where $1 \le i \le 3$ and $0 \leq  k < \ell \le 2$.
\end{enumerate}
Furthermore, $C_{1,2}^{(i)} \sim 3\beta_3$
and $C_{0,1}^{(i)} \sim C_{0,2}^{(i)} \sim \beta_3$ for every $i$.
\end{lemma}
\begin{proof}
It remains to prove the last sentence. Identify $M_3(P_i)$ with
the punctual Hilbert scheme $\Hilb^3(\C^2, 0)$.
By (\ref{tilde_sigma_n:def}), $C_{1,2}^{(i)} = \tilde \sigma_3$.
It follows from Lemma~\ref{tilde_sigma_n:homology} that
$C_{1,2}^{(i)} \sim 3\beta_3$. Similarly, we see from
(\ref{sigma_n:def}) and Lemma~\ref{sigma_n:homology} that
$C_{0,1}^{(i)} \sim C_{0,2}^{(i)} \sim \beta_3$.
\end{proof}

Next, we compute the equivariant first Chern classes of
the restrictions of the tautological bundles (\ref{E_i})
to the $T$-fixed points in $(\PP^2)^{[3]}$.
Let $w_i = c_1(\lambda_i)$ and $z_i = c_1(\mu_i)$ in
the equivariant Chow group $A^T_*(pt)$. If we put $(w_0, z_0) = (w,z)$,
then $(w_1,z_1) = (-w,-w + z)$ and $(w_2,z_2) = (-z, -z+w)$.

\begin{lemma} \label{lem.fixedpoints}
Let $g_{0}=0$, $g_{1}= w$, and $g_{2}= z$.
There are $T$-linearizations on $\cE_0$ and $\cE_1$ such that
$c_1(\cE_0|_{R_{i,j}^{(1)}}) = z_i$,
$c_1(\cE_0|_{R_{i,j}^{(2)}}) = w_i$,
$c_1(\cE_0|_{Q_{i,0}}) = z_i + w_i$,
$c_1(\cE_0|_{Q_{i,1}}) = 3z_i$,
$c_1(\cE_0|_{Q_{i,2}}) = 3w_i$ and
$c_1(\cE_1|_{R_{i,j}^{(1)}}) = 2g_{i} + g_{j} + z_i$,
$c_1(\cE_1|_{R_{i,j}^{(2)}}) =  2g_{i} + g_{j} + w_i$,
$c_1(\cE_1|_{Q_{i,0}}) = 3g_{i} + z_i + w_i$,
$c_1(\cE_1|_{Q_{i,1}}) = 3g_{i} + 3z_i$,
$c_1(\cE_1|_{Q_{i,2}}) = 3g_{i} + 3w_i$.
\end{lemma}
\begin{proof}
The proofs of these conclusions are similar. For instance,
let us prove $c_1(\cE_1|_{R_{i,j}^{(2)}}) =  2g_{i} + g_{j} + w_i$.
Note that the fiber $\cE_1|_{R_{i,j}^{(2)}}$ is canonically identified
with $\mathcal O_X(1) \otimes \mathcal O_X/I_{R_{i,j}^{(2)}}$.
Since $R_{i,j}^{(2)} = \xi_{i, 2} + P_j$,
$\cE_1|_{R_{i,j}^{(2)}}$ is canonically identified with
$\left ( \mathcal O_X(1) \otimes \mathcal O_X/I_{\xi_{i, 2}} \right )
\oplus \left ( \mathcal O_X(1) \otimes \mathcal O_X/I_{P_j} \right )$.
Therefore,
\begin{eqnarray*}
c_1(\cE_1|_{R_{i,j}^{(2)}}) = 2c_1(\mathcal O_X(1)|_{P_i})
+c_1(\mathcal O_X/I_{\xi_{i, 2}}) +c_1(\mathcal O_X(1)|_{P_j}).
\end{eqnarray*}
Since $\mathcal O_X(1)|_{P_i} \cong (\C \oplus \C)/(\C P_i)$,
we have $c_1(\mathcal O_X(1)|_{P_i}) = g_{i}$.
Using $c_1(\mathcal O_X/I_{\xi_{i, 2}}) = c_1(\lambda_i) = w_i$,
we conclude that $c_1(\cE_1|_{R_{i,j}^{(2)}}) =  2g_{i} + g_{j} + w_i$.
\end{proof}

\subsection{The Euler characteristic for a covering}
\label{subsect_euler_covering}
\par
$\,$

An important step in computing the virtual Euler class
of the $T$-fixed locus $\M_{0,2}((\PP^2)^{[3]}, d\beta_3)^T$
is to compute (as a representation) $\chi (f^*T_{(\PP^2)^{[3]}})$
where $f: \PP^1 \to (\PP^2)^{[3]}$ is a degree-$d$ morphism such that
the image is one of the 15 $T$-invariant curves in Lemma~\ref{15_curve}
and $f$ is totally ramified at the two $T$-fixed points in $f(\PP^1)$.

\subsubsection{\bf Degree-$d$ coverings of $C_{k,\ell}^{(i)}$}
\par
$\,$

Observe that if $\PP^1 \to (\PP^2)^{[3]}$ is a degree-$d$
$T$-equivariant morphism with image $C_{k,\ell}^{(i)}$,
then the characters of $T$-action on $\PP^1$ are (using multiplicative
notation) $\alpha^{1/d}, \beta^{1/d}$ where $\alpha, \beta$ are
the characters of the $T$-action on the image curve $C_{k,\ell}^{(i)}$.
Let $S_{i, k}$ and $S_{i, \ell}$ be the two fixed points of the action
on $\PP^1$ denoted so that the image of $S_{i, k}$ is
$Q_{i, k}$ and the image of $S_{i, \ell}$ is $Q_{i, \ell}$.
If $V$ is a $T$-equivariant vector bundle on $\PP^1$,
then the localization theorem for equivariant $K$-theory says that
\begin{equation} \label{eq.locformula}
\chi(V) = \frac{V|_{S_{i, k}}}{1 - T^*_{\PP^1}|_{S_{i, k}}} +
\frac{V|_{S_{i, \ell}}}{1 - T^*_{\PP^1}|_{S_{i, \ell}}}
\end{equation}
where $T^*_{\PP^1}$ is the cotangent bundle of $\PP^1$.
Since $T^*_{\PP^1}|_{S_{i, k}} \cong T^*_{C_{k, \ell}^{(i)}}|_{Q_{i, k}}$,
we can use formulas \eqref{eq.tq0}, \eqref{eq.tq1}, \eqref{eq.tq2}
to determine $\chi (f^*T_{(\PP^2)^{[3]}})$.

First of all, let $f(\PP^1) = C_{0,1}^{(i)}$.
The curve $C_{0,1}^{(i)}$ is a component of the fixed locus of
$\ker(\lambda_i\mu_i^{-2})$.
Thus, reading off (\ref{eq.tq0}) and (\ref{eq.tq1}), we see that
$T_{{C_{0, 1}^{(i)}}}|_{Q_{i, 0}}= \lambda_i\mu_i^{-2}$
and $T_{{C_{0, 1}^{(i)}}}|_{Q_{i, 1}}= \lambda_i^{-1}\mu_i^{2}$. Thus
$T_{\PP^1}|_{S_{i, 0}} =\gamma_i \theta_i^{-2}$ and
$T_{\PP^1}|_{S_{i, 1}} =\gamma_i^{-1}\theta_i^{2}$ where
$\gamma_i^d = \lambda_i$ and $\theta_i^{d} = \mu_i$.
Substituting \eqref{eq.tq0} and \eqref{eq.tq1} into
the localization formula \eqref{eq.locformula} yields
\begin{eqnarray*}
 \chi (f^*T_{(\PP^2)^{[3]}})
&= &\frac{\lambda_i\mu_i^{-2} + \mu_i^{-1} + \lambda_i^{-1}
 + \lambda_i^{-2}\mu_i +\lambda_i^{-1}
 + \mu_i^{-1}}{1- \gamma_i^{-1}\theta_i^2}  \\
&+ &\frac{\lambda_i^{-1}\mu_i^{2} + \lambda_i^{-1}\mu_i + \lambda_i^{-1}
 + \mu_i^{-3} +\mu_i^{-2} + \mu_i^{-1}}{1 - \gamma_i\theta_i^{-2}}.
\end{eqnarray*}
Since $1/(1- \gamma_i^{-1}\theta_i^2)
= -\gamma_i\theta_i^{-2}/(1 - \gamma_i\theta_i^{-2})$,
the right hand side can be rewritten as
\begin{eqnarray*}
&\displaystyle{\frac{1}{1-\gamma_i\theta_i^{-2}} \big [
(\lambda_i^{-1}\mu_i^2 + \lambda_i^{-1}\mu_i
+ \lambda_i^{-1}+ \mu_i^{-3}+ \mu_i^{-2} + \mu_i^{-1})}
- \gamma_i\theta_i^{-2} \big (
(\lambda_i^2\mu_i^{-4})(\lambda_i^{-1}\mu_i^2) \\
&+ (\lambda_i\mu_i^{-2})\lambda_i^{-1}\mu_i + \lambda_i^{-1}
+ (\lambda_i^{-2}\mu_i^4) \mu_i^{-3} + (\lambda_i^{-1}\mu_i^2)
\mu_i^{-2} + \mu_i^{-1} \big ) \big ].
\end{eqnarray*}
Using $\lambda_i= \gamma_i^d$ and $\mu_i= \theta_i^d$,
we conclude that $\chi(f^*T_{(\PP^2)^{[3]}})$ is equal to
\begin{eqnarray*}
&\lambda_i^{-1}\mu_i^2 \sum\limits_{m=0}^{2d} (\gamma_i\theta_i^{-2})^m
   + \lambda_i^{-1}\mu_i \sum\limits_{m=0}^{d}
   (\gamma_i\theta_i^{-2})^m +\lambda_i^{-1}&\\
&-\mu_i^{-3}(\gamma_i\theta_i^{-2})^{-2d +1}
  \sum\limits_{m=0}^{2d-2} (\gamma_i\theta_i^{-2})^m
  -\mu_i^{-2}(\gamma_i\theta_i^{-2})^{-d +1}
  \sum\limits_{m=0}^{d-2} (\gamma_i\theta_i^{-2})^m + \mu_i^{-1}.&
\end{eqnarray*}
To simplify this further, set $\Theta_{0,1}^{(i)} = \sum_{m= 1}^{d-1}
(\gamma_i\theta_i^{-2})^m= \sum_{m= 1}^{d-1} (\lambda_i\mu_i^{-2})^{m/d}$
(with the understanding that $\Theta_{0,1}^{(i)} = 0$ when $d=1$).
Then we see that $\chi(f^*T_{(\PP^2)^{[3]}})$ equals
%\begin{eqnarray*}
%\chi(f^*T_{(\PP^2)^{[3]}})
%& = &(\lambda_i^{-1}\mu_i^2 + \lambda_i^{-1}\mu_i^2\Theta_{0,1}^{(i)} + 1 +
%\Theta_{0,1}^{(i)} +
%\lambda_i\mu_i^{-2})\\
%& + & (\lambda_i^{-1}\mu_i + \lambda_i^{-1}\mu_i\Theta_{0,1}^{(i)}
%+ \mu_i^{-1}) + \lambda_i^{-1} \\
%& - & (\lambda_i^{-2}\mu_i\Theta_{0,1}^{(i)} + \lambda_i^{-1}\mu_i^{-1}
%+ \lambda_i^{-1}\mu_i^{-1}\Theta{0,1}^{(i)})-
%\lambda_i^{-1}\Theta_{0,1}^{(i)}+ \mu_i^{-1},
%\end{eqnarray*}
\begin{eqnarray} \label{eq.chic01}
&(1 + \lambda_i^{-1}\mu_i^2 +\lambda_i\mu_i^{-2} + \lambda_i^{-1}\mu_i
+ \mu_i^{-1} + \lambda_i^{-1} + \mu_i^{-1}
- \lambda_i^{-1}\mu_i^{-1})&  \nonumber \\
& + (\lambda_i^{-1}\mu_i^2 + 1 + \lambda_i^{-1}\mu_i - \lambda_i^{-2}\mu_i
- \lambda_i^{-1}\mu_i^{-1} - \lambda_i^{-1})\Theta_{0,1}^{(i)}.&
\end{eqnarray}

By symmetry, if $f(\PP^1) = C_{0,2}^{(i)}$,
then $\chi(f^*T_{(\PP^2)^{[3]}})$ is equal to
\begin{eqnarray}   \label{eq.chic02}
&(1 + \mu_i^{-1}\lambda_i^2 +\mu_i\lambda_i^{-2} + \mu_i^{-1}\lambda_i
+ \lambda_i^{-1} + \mu_i^{-1} + \lambda_i^{-1}
- \mu_i^{-1}\lambda_i^{-1})&    \nonumber    \\
& + (\mu_i^{-1}\lambda_i^2 + 1 + \mu_i^{-1}\lambda_i - \mu_i^{-2}\lambda_i
- \mu_i^{-1}\lambda_i^{-1} - \mu_i^{-1})\Theta_{0,2}^{(i)}&
\end{eqnarray}
where $\Theta_{0,2}^{(i)} = \sum_{m=1}^{d-1} (\mu_i\lambda_i^{-2})^{m/d}$,
and as above $\Theta_{0,2}^{(i)} = 0$ if $d = 1$.

Next, let $f(\PP^1) = C_{1,2}^{(i)}$. Then
$T_{{C_{1,2}^{(i)}}}|_{Q_{i, 1}}= \lambda_i^{-1}\mu_i$
and $T_{{C_{1,2}^{(i)}}}|_{Q_{i, 2}}= \lambda_i\mu_i^{-1}$.
Thus $T_{\PP^1}|_{S_{i, 1}} =\gamma_i^{-1} \theta_i$ and
$T_{\PP^1}|_{S_{i, 2}} =\gamma_i\theta_i^{-1}$.
By (\ref{eq.locformula}), (\ref{eq.tq1}) and (\ref{eq.tq2}),
$\chi(f^*T_{(\PP^2)^{[3]}})$ equals
\begin{eqnarray*}
&\displaystyle{\frac{1}{1 - \gamma_i\theta_i^{-1}}[
(\lambda_i^{-1}\mu_i^2 + \lambda_i^{-1}\mu_i + \lambda_i^{-1} +
\mu_i^{-3} + \mu_i^{-2} + \mu_i^{-1})}  \\
&- \gamma_i\theta_i^{-1}(\lambda_i^2 \mu_i^{-1} + \lambda_i\mu_i^{-1}
+ \mu_i^{-1} + \lambda_i^{-3} + \lambda_i^{-2} + \lambda_i^{-1})]
\end{eqnarray*}
As above, the numerator is divisible by $(1 - \gamma_i\theta_i^{-1})$,
and $\chi(f^*T_{(\PP^2)^{[3]}})$ is equal to
\begin{eqnarray*}
\lambda_i^{-1} \sum\limits_{s=0}^2 \mu_i^s
 \sum\limits_{m=0}^{(s+1)d} (\gamma_i\theta_i^{-1})^m
 - \sum\limits_{s=1}^3 \lambda_i^{-s}
 \sum\limits_{m=1}^{sd-1} (\gamma_i\theta_i^{-1})^m.
\end{eqnarray*}
Let $\Theta_{1,2}^{(i)} = \sum_{m= 1}^{d-1} (\lambda_i\mu_i^{-1})^{m/d}$
with $\Theta_{1,2}^{(i)} = 0$ when $d=1$.
Then $\chi(f^*T_{(\PP^2)^{[3]}})$ equals
\begin{eqnarray*}
&\lambda_i^{-1}(1 + \Theta_{1,2}^{(i)}) + \mu_i^{-1}
 + \lambda_i^{-1}\mu_i(1 + \theta_{1,2}^{(i)}) + (1 + \theta_{1,2}^{(i)})
 + \lambda_i\mu_i^{-1} &  \\
&\lambda_i^{-1}\mu_i^2(1 + \Theta_{1,2}^{(i)})+ \mu_i(1 + \Theta_{1,2}^{(i)})
 + \lambda_i(1 + \Theta_{1,2}^{(i)}) + \lambda_i^2\mu_i^{-1}&  \\
&- \lambda_i^{-1}\Theta_{1,2}^{(i)} - (\lambda_i^{-2}\Theta_{1,2}^{(i)}
 + \lambda_i^{-1}\mu_i^{-1}(1 + \Theta_{1,2}^{(i)}))  &  \\
&- (\lambda_i^{-3}\Theta_{1,2} + \lambda_i^{-2}\mu_i^{-1}(1 + \Theta_{1,2}^{(i)})
 + \lambda_i^{-1}\mu_i^{-2}(1 + \Theta_{1,2}^{(i)})).  &
\end{eqnarray*}
Rearranging the terms, we conclude that $\chi(f^*T_{(\PP^2)^{[3]}})$
is equal to
\begin{eqnarray}   \label{eq.chic12}
&(\lambda_i^{-1} + \mu_i^{-1}+ \lambda_i^{-1}\mu_i+1 + \lambda_i\mu_i^{-1}
  + \lambda_i^{-1}\mu_i^2 + \mu_i + \lambda_i + \lambda_i^2\mu_i^{-1} & \\
&- \lambda_i^{-1}\mu_i^{-1} - \lambda_i^{-2}\mu_i^{-1}
  -\lambda_i^{-1}\mu_i^{-2})&   \nonumber \\
&+ ( 1 + \lambda_i^{-1}\mu_i^2 + \mu_i +\lambda_i
  + \lambda_i^{-1}\mu_i - \lambda_i^{-2} - \lambda_i^{-1}\mu_i^{-1}
  - \lambda_i^{-3} - \lambda_i^{-2}\mu_i^{-1}- \lambda_i^{-1}\mu_i^{-2})
  \Theta_{1,2}^{(i)}.  &  \nonumber
\end{eqnarray}

\subsubsection{\bf Degree-$d$ coverings of $C_{i,j}$}
\par
$\,$

Consider maps $f:\PP^1 \to (\PP^2)^{[3]}$ which
are degree-$d$ and have image $C_{i,j}$.
To compute $\chi(f^*T_{(\PP^2)^{[3]}})$,
we recall from subsection \ref{subsect_curve} that
the $T$-fixed points on $C_{i,j}$ are
$R_{i,j}^{(1)}$ and $R_{i,j}^{(2)}$.
Using the results in \cite{E-S}, we have the following decompositions
for the tangent spaces of $(\PP^2)^{[3]}$ at $R_{i,j}^{(1)}$ and
$R_{i,j}^{(2)}$ as representations of $T$:
\begin{eqnarray}
T_{R_{i,j}^{(1)}} &=& \lambda_i^{-1}\mu_i + \lambda_i^{-1} + \mu_i^{-2}
+ \mu_i^{-1} + \lambda_j^{-1}+ \mu_j^{-1}, \label{eq.tr1}  \\
T_{R_{i,j}^{(2)}} &=& \lambda_i^{-2} + \lambda_i^{-1}+ \lambda_i\mu_i^{-1}
+ \mu_i^{-1} + \lambda_j^{-1} + \mu_j^{-1}. \label{eq.tr2}
\end{eqnarray}
%(Note the last two terms are just the tangent space to $X$
%at $P_j$. We take the reciprocal of the weights on the cordinates
%to get the weights on the tangent space. This also follows
%from plugging into the formula of \cite{E-S} when $d = 1$.)
Also, $T_{C_{i,j}}|_{R_{i,j}^{(1)}} = \lambda_i^{-1}\mu_i$ and
$T_{C_{i,j}}|_{R_{i,j}^{(2)}} = \lambda_i\mu_i^{-1}$.
By (\ref{eq.locformula}), $\chi (f^*T_{(\PP^2)^{[3]}})$ equals
\begin{eqnarray*}
\quad \frac{\lambda_i^{-1}\mu_i + \lambda_i^{-1} + \mu_i^{-2}
+ \mu_i^{-1} + \lambda_j^{-1}+ \mu_j^{-1}}{1- \gamma_i\theta_i^{-1}}
+ \frac{\lambda_i^{-2} + \lambda_i^{-1}+ \lambda_i\mu_i^{-1}
+ \mu_i^{-1} + \lambda_j^{-1} + \mu_j^{-1}}{1 - \gamma_i^{-1}\theta_i }.
\end{eqnarray*}
So we obtain the following formula for $\chi (f^*T_{(\PP^2)^{[3]}})$:
\begin{eqnarray} \label{eq.chicij}
&(1 + \lambda_i^{-1}\mu_i +\lambda_i\mu_i^{-1} + \lambda_i^{-1}
+ \mu_i^{-1} + \lambda_j^{-1} + \mu_j^{-1}
- \lambda_i^{-1}\mu_i^{-1})&  \nonumber \\
& + (1 + \lambda_i^{-1}\mu_i - \lambda_i^{-2} - \lambda_i^{-1}\mu_i^{-1}
)\Theta_{1,2}^{(i)}.&
\end{eqnarray}

\subsection{$T$-invariant stable maps, stable graphs and localizations}
\label{subsect_graph}
\par
$\,$

Let $X = \PP^2$. Note that if $[f: (C; p_1, p_2) \to X^{[3]}] \in
\overline {\frak M}_{0, 2}(X^{[3]}, d\beta_3)$ is $T$-invariant
and if $\mathbb P^1$ is an irreducible component of $C$ with
nonconstant $f|_{\mathbb P^1}$, then $f({\mathbb P^1})$ is one of
the $15$ $T$-invariant curves in Lemma~\ref{15_curve}.
The restriction $f|_{\mathbb P^1}$ is ramified at
exactly two points with ramification index $\text{deg}(f|_{\mathbb P^1})$.
Since $f|_{\mathbb P^1}$ is ramified at every special point,
$\mathbb P^1$ contains at most two special points. Moreover, $f$ maps
the contracted components and the special points
(i.e., marked points, nodal points and ramification points)
of $C$ into the $T$-fixed point set $(X^{[3]})^T$.

Following the book \cite{C-K}, to each $T$-invariant stable map
$[f: (C; p_1, p_2) \to X^{[3]}] \in
\overline {\frak M}_{0, 2}(X^{[3]}, d\beta_3)$,
we can associate a marked graph $\Gamma$
called {\it a stable graph of genus-$0$}. The graph $\Gamma$ has
one vertex for each connected component of $f^{-1}((X^{[3]})^T)$.
It has one edge $e$ for each non-contracted component $C_e \simeq \mathbb P^1$,
whose two vertices correspond to the connected components of
$f^{-1}((X^{[3]})^T)$ containing the two ramification points
in the component $C_e$. The edge $e$ is marked with the degree
$d_e \,\, {\buildrel\text{def}\over=} \,\, \text{deg}(f|_{C_e})$.
Note that the morphism $f$ defines a labeling map $\mathfrak L$
from the vertices of $\Gamma$ to $(X^{[3]})^T$.
Finally, a vertex is marked with
$\{1 \}$ (respectively, $\{ 2 \}$, or $\{1, 2 \}$) if the connected component
of $f^{-1}((X^{[3]})^T)$ corresponding to the vertex contains
the marked point $p_1$ (respectively, $p_2$, or both $p_1$ and $p_2$).

To a stable graph $\Gamma$, we introduce the following notation
(cf. \cite{C-K}). Recall that a flag $F$ is a pair $(v, e)$ consisting of
an edge $e$ and a vertex $v$ of $e$.
For a flag $F=(v,e)$, define $i(F) = {\mathfrak
L}(v)$. Let $S(v)$ be the number of markings of $v$, and $val(v)$ be
the valance of $v$ (i.e., the number of edges $e$ such that $v$ is a
vertex of $e$). Let $n(F) = n(v)= val(v) + S(v)$.  If $val(v) =1$, let
$F(v)$ be the single flag containing $v$; if $val(v) = 2$, let
$F_1(v)$ and $F_2(v)$ denote the two flags containing $v$.

Now the connected components of
$\overline{\frak M}_{0,2}(X^{[3]}, d\beta_3)^T$
are enumerated by stable graphs corresponding to stable maps
whose images are unions of the $15$ $T$-invariant curves
in Lemma~\ref{15_curve} and whose contracted components and
special points are mapped into $(X^{[3]})^T$.
We use $\Gamma$ to denote these stable graphs,
and use ${\fM}_\Gamma$ to denote the corresponding connected components
of $\overline{\frak M}_{0,2}(X^{[3]}, d\beta_3)^T$.
If $\Gamma$ is a stable graph,
let $M_\Gamma = \prod_{n(v) \geq 3} \overline{M}_{0,n(v)}$
where $\overline{M}_{0,n(v)}$ is the (fine) moduli space of
$n(v)$-pointed stable rational curves. As discussed in \cite{C-K},
there is a finite map $M_\Gamma \to {\mathfrak M}_\Gamma$ such that
${\mathfrak M}_\Gamma = M_\Gamma/{\bf A}_\Gamma$ where
${\bf A}_\Gamma$ fits in the exact sequence
\begin{eqnarray*}
0 \to \prod_{e} {\mathbb Z}/d_e {\mathbb Z} \to {\bf A}_\Gamma \to
{\text{\rm Aut}} (\Gamma) \to 0.
\end{eqnarray*}

Since a stable curve is connected, we see from the description of
the $T$-invariant curves in Lemma~\ref{15_curve} that
a summation over all the stable graphs $\Gamma$ breaks up as
\begin{eqnarray}  \label{break}
\sum_\Gamma \,\, = \,\, \sum_{1\leq i\neq j \leq 3} \,\,
\sum_{\Gamma \in {\mathcal S}_{d,i,j}} \,\,
+\,\, \sum_{i=1}^{3} \,\,\sum_{\Gamma \in \mathcal T_{d,i}}
\end{eqnarray}
where ${\mathcal S}_{d,i,j}$ is the set of all stable graphs $\Gamma$
such that $f(C) = C_{i,j}$ for every $[f: (C; p_1, p_2) \to X^{[3]}]
\in {\fM}_\Gamma$, and $\mathcal T_{d,i}$ is the set of
all stable graphs $\Gamma$ such that
$f(C) \subset C_{0,1}^{(i)} \cup C_{0,2}^{(i)} \cup C_{1,2}^{(i)}$
for every $[f: (C; p_1, p_2) \to X^{[3]}] \in {\fM}_\Gamma$.

Our goal of this section is to study $\langle \pd(\frak a_{-3}(\ell)\vac),
\pd(\frak a_{-3}(X)\vac) \rangle_{0, d}$.
To apply the localization formula more effectively,
we rewrite this $2$-point invariant by using the Chern classes of
tautological bundles over $X^{[3]} = (\PP^2)^{[3]}$ defined in (\ref{E_i}).
Let
\begin{eqnarray*}
A = (c_1(\mathcal E_1)- c_1(\mathcal E_0)) c_1(\mathcal E_0)^2
\quad \text{\rm and} \quad B = c_1(\mathcal E_0)^2.
\end{eqnarray*}
Intersecting (\ref{a_3_X}) with
$D_\ell = c_1(\mathcal E_1)- c_1(\mathcal E_0)$, we see that
$A$ is equal to
\begin{eqnarray*}
&3\frak a_{-3}(\ell)\vac - 3
   \frak a_{-1}(X)\frak a_{-1}(\ell)\frak a_{-1}(x)\vac & \\
&- {1 \over 2}\frak a_{-1}(\ell)^3\vac
   + 3\frak a_{-1}(X)\frak a_{-2}(x)\vac
   + {3 \over 2}\frak a_{-2}(\ell)\frak a_{-1}(\ell)\vac.&
\end{eqnarray*}
By Lemma~\ref{2_zero}, Lemma~\ref{2_zero_K3} and Lemma~\ref{2_obs}~(i),
we obtain
\begin{eqnarray}  \label{tilde_A_to_A}
\langle A, B \rangle_{0, d} = 3 \,\, \langle \pd(\frak a_{-3}(\ell)\vac),
\pd(\frak a_{-3}(X)\vac) \rangle_{0, d}
\end{eqnarray}
where for notational simplicity, we make no distinction between
the algebraic cycles $A, B$ and their corresponding cohomology classes.

By the virtual localization formula of \cite{G-P}, we have
\begin{eqnarray}\label{eq.gwvirtloc}
   \langle A, B \rangle_{0, d}
   =\int_{[\overline{\frak M}_{0, 2}(X^{[3]},d\beta_3)]^{\rm vir}}
    ev_2^*(A \otimes B)
=\sum_{\Gamma} \frac{1}{|{\bf A}_\Gamma|}\int_{[M_\Gamma]^{\rm vir}}
   \frac{(A \otimes B)_\Gamma}{e(N_\Gamma^{\rm vir})}.
\end{eqnarray}
%Here the summation $\sum_{\Gamma}$ is indexed by the connected components
%$\fM_\Gamma$ of the fixed locus
%$\overline {\frak M}_{0, 2}(X^{[3]}, d\beta_3)^T$, and
Here $[M_\Gamma]^{\rm vir}$ is the pullback of
$[\fM_\Gamma]^{\rm vir}$ to $M_\Gamma$ via the finite map
$M_\Gamma \to \fM_\Gamma$.
Likewise, $(A \otimes B)_\Gamma$ is the pullback of
$ev_2^*(A \otimes B)|_{\fM_\Gamma}$ to $M_\Gamma$,
and $e(N_\Gamma^{\rm vir})$ is
the pullback of the Euler class of the moving part $N_\Gamma^{\rm vir}$ of
the tangent-obstruction complex.

Let $\Gamma$ be a stable graph such that the labeling $\mathfrak L$
maps the marked vertices of $\Gamma$ to the same point in $(X^{[3]})^T$.
Then we have $(A \otimes B)_\Gamma = (1_X \otimes AB)_\Gamma$
where $1_X \in H^0(X)$ is the fundamental cohomology class.
By the fundamental class axiom, $\langle 1_X, AB \rangle_{0, d} =0$.
Thus in view of (\ref{eq.gwvirtloc}) and (\ref{break}), we obtain
\begin{eqnarray}  \label{1_X_prime}
& &\langle A, B \rangle_{0, d} = \langle A, B \rangle_{0, d}
   - \langle 1_X, AB \rangle_{0, d}    \nonumber \\
&=&{\sum_{\Gamma}} \int_{[M_\Gamma]^{\rm vir}}
   \frac{(A \otimes B)_\Gamma -
(1_X \otimes AB)_\Gamma}{|{\bf A}_{\Gamma}| \,\, e(N_\Gamma^{\rm vir})}
   = \,\, \sum_{1\leq i\neq j \leq 3} \,\,
   \sum_{\Gamma \in {\mathcal S}_{d,i,j}'} \,\, +\,\,
   \sum_{i=1}^{3} \,\,\sum_{\Gamma \in \mathcal T_{d,i}'} \quad \qquad
\end{eqnarray}
where the three prime signs indicate that we only sum over stable graphs
$\Gamma$ such that the two
marked vertices of $\Gamma$ have distinct labels in $(X^{[3]})^T$.
In other words, putting $S_{d,i,j}' = \sum_{\Gamma \in {\mathcal S}_{d,i,j}'}$
and $T_{d,i}' = \sum_{\Gamma \in \mathcal T_{d,i}'}$, we have
\begin{eqnarray} \label{prime}
\langle A, B \rangle_{0, d} \,\,= \,\, \sum_{1\leq i\neq j \leq 3} \,\,
S_{d,i,j}'  \,\, +\,\, \sum_{i=1}^{3} \,\, T_{d,i}'.
\end{eqnarray}

\subsection{Computation of $S_{d,i,j}'$} \label{subsect_Sdij'}
\par
$\,$

Let ${\mathcal S}_{d,i,j}'' = {\mathcal S}_{d,i,j}'/\sim$
where $\Gamma_1 \sim \Gamma_2$ if $\Gamma_1$ and $\Gamma_2$ are identical
except that the vertex which is marked with $\{1\}$
(respectively, with $\{2\}$) in $\Gamma_1$ is marked with
$\{2\}$ (respectively, with $\{1\}$) in $\Gamma_2$.
%one can be obtained from the other
%by switching the two marked vertices.
Then each graph $\Gamma$ in ${\mathcal S}_{d,i,j}''$ gives rise to
two graphs $\Gamma_1, \Gamma_2$ in ${\mathcal S}_{d,i,j}'$.
However, there is no ambiguity to define
\begin{eqnarray} \label{edij}
e_{d,i,j} = \sum_{\Gamma \in {\mathcal S}_{d,i,j}''}
\int_{[M_{\Gamma_1}]^{\rm vir}} \frac{1}{|{\bf A}_{\Gamma_1}|
\,\, e(N_{\Gamma_1}^{\rm vir})}.
\end{eqnarray}
By the definition of ${\mathcal S}_{d,i,j}$,
$f(C) = C_{i,j}$ for every stable map $[f: (C; p_1, p_2) \to X^{[3]}]$
in ${\fM}_{\Gamma_1}$ or ${\fM}_{\Gamma_2}$. Recall that $R_{i,j}^{(1)}$
and $R_{i,j}^{(2)}$ are the two $T$-fixed points in $C_{i,j}$. So
\begin{eqnarray*}
& &\int_{[M_{\Gamma_1}]^{\rm vir}} \frac{(A \otimes B)_{\Gamma_1} -
   (1_X \otimes AB)_{\Gamma_1}}{|{\bf A}_{\Gamma_1}| \,\,
   e(N_{\Gamma_1}^{\rm vir})}
   + \int_{[M_{\Gamma_2}]^{\rm vir}} \frac{(A \otimes B)_{\Gamma_2} -
   (1_X \otimes AB)_{\Gamma_2}}{|{\bf A}_{\Gamma_2}| \,\,
   e(N_{\Gamma_2}^{\rm vir})} \nonumber \\
&=&-(A|_{R_{i,j}^{(1)}} - A|_{R_{i,j}^{(2)}})
   (B|_{R_{i,j}^{(1)}} - B|_{R_{i,j}^{(2)}})
   \cdot \int_{[M_{\Gamma_1}]^{\rm vir}} \frac{1}{|{\bf A}_{\Gamma_1}|
   \,\, e(N_{\Gamma_1}^{\rm vir})}.
\end{eqnarray*}
Combining this with Lemma \ref{lem.fixedpoints} and (\ref{edij}),
we conclude that
\begin{eqnarray} \label{Gamma_12}
S_{d,i,j}'= -(2g_{i} + g_{j})(w_i^2 - z_i^2)^2 \,\, e_{d,i,j}.
\end{eqnarray}

To compute $e_{d,i,j}$, we calculate the contribution from
a graph $\Gamma_1$ by considering the restriction of
the tangent-obstruction complex
on $\overline{\frak M}_{0,2}(X^{[3]}, d\beta_3)$ to $\fM_{\Gamma_1}$.
Following \cite{G-P}, the fibers of its cohomology sheaves,
${\mathcal T}^1$ and ${\mathcal T}^2$,
at a point associated to a stable map
$[f: (C; p_1, p_2) \to X^{[3]}]$ fit into the exact sequence
$$\begin{array}{ccccccccc}
0 & \to  & \Ext^0(\Omega_C(p_1 + p_2),{\mathcal O}_C)& \to &
H^0(C,f^*T_{X^{[3]}}) & \to & {\mathcal T}^1 \\
& \to & \Ext^1(\Omega_C(p_1 + p_2),{\mathcal O}_C) & \to &
H^1(C,f^*T_{X^{[3]}}) & \to & {\mathcal T}^2  & \to & 0.
\end{array}$$
To obtain the contribution of the moving parts of each term in
the sequence, we use an analysis similar to that carried out for $\PP^r$
in \cite{G-P}. As was the case for $\PP^r$,
the fixed part ${\mathcal T}^{2, \rm f}$ vanishes.
So the fixed stack is smooth with tangent bundle ${\mathcal T}^{1, \rm f}$.
In particular $[{\fM}_{\Gamma_1}]^{\rm vir} = [\fM_{\Gamma_1}]$.
%The only difference is that, $X^{[3]}$ is not convex,
%so $H^1(C,f^*|_{C_e}TX^{[3]})\neq 0$
%where $C_e$ is component corresponding to an edge (i.e.,
%a component which is not contracted).
As a result, denoting the contributions from the edges, vertices and flags
of the graph $\Gamma_1$ by $e_{\Gamma_1}^{\rm e},
e_{\Gamma_1}^{\rm v}, e_{\Gamma_1}^{\rm F}$ respectively, we obtain
\begin{eqnarray}  \label{e_v_F}
e(N_{\Gamma_1}^{\rm vir}) = e_{\Gamma_1}^{\rm e} \cdot
e_{\Gamma_1}^{\rm v} \cdot e_{\Gamma_1}^{\rm F}.
\end{eqnarray}

First of all, we have $e_{\Gamma_1}^{\rm e} =
\prod_e e(\chi(((f|_{C_e})^*T_{X^{[3]}})^{\rm m}))$
where $((f|_{C_e})^*T_{X^{[3]}})^{\rm m}$ denotes the moving part
in $(f|_{C_e})^*T_{X^{[3]}}$.
%This can be calculated from formula (\ref{eq.chicij}) by noting
%that all terms except $1$ contribute to the moving part.
It follows from (\ref{eq.chicij}) that
\begin{eqnarray}  \label{e_v_F.e}
e_{\Gamma_1}^{\rm e} = \prod_e
\frac{(-1)^{d_e -1}((d_e-1)!)^2w_i w_j z_i z_j (w_i - z_i)^2}
{(w_i + z_i) P(1 + \frac{2d_ew_i}{-w_i + z_i},d_e -1)
P(1 - \frac{d_e(w_i + z_i)}{w_i - z_i},d_e - 1)}
\end{eqnarray}
where $P(a,n)$ denotes the polynomial $a(a +1) \ldots (a + n-1)$.

Now the contributions of vertices and flags are given by
\begin{eqnarray}
e_{\Gamma_1}^{\rm v}
&=&\prod_v e(T_{{\frak L}(v)}) \cdot
   \prod_{val(v) = n(v) =2}(\omega_{F_1(v)}+ \omega_{F_2(v)}) \cdot
   \prod_{val(v) = n(v) = 1} \omega_{F(v)}^{-1}  \label{e_v_F.v}\\
e_{\Gamma_1}^{\rm F}
&=&\prod_{n(F) \geq 3}(\omega_F - e_F)
   \cdot \prod_F e(T_{i(F)})^{-1}   \label{e_v_F.F}
\end{eqnarray}
where for a flag $F = (v, e)$,
we put $\omega_F = e(T_{i(F)}C_{i,j})/d_e$,
and define $e_F$ to be the first Chern class of the bundle on
$M_\Gamma$ whose fiber is the cotangent space of the component
associated to $v$ at the point corresponding to the flag $F$
(c.f. \cite[p.285]{C-K}).
Note that $T_{i(F)} = T_{{\frak L}(v)}$ has been computed in
(\ref{eq.tr1}) and (\ref{eq.tr2}).
Thus, $\omega_F= (-w_i + z_i)/d_e$ if $i(F) = R_{i,j}^{(1)}$,
and $\omega_F = (w_i - z_i)/d_e$  if $i(F) = R_{i,j}^{(2)}$.

%$e(T_{R_{i,j}^{(1)}}) = 2w_iw_j(w_i - z_i)z_i^2z_j$
%and $e(T_{R_{i,j}^{(2)}}) = -2w_i^2w_j(w_i - z_i)z_iz_j$.

%
%
%
%
%
%
%
%
%
%
\subsection{Computation of $T_{d,i}'$} \label{subsect_Tdi'}
\par
$\,$

Recall from (\ref{1_X_prime}) and (\ref{break}) that $\mathcal T_{d,i}'$
is the set of all stable graphs $\Gamma$ such that
$f(C) \subset C_{0,1}^{(i)} \cup C_{0,2}^{(i)} \cup C_{1,2}^{(i)}$
for every $[f: (C; p_1, p_2) \to X^{[3]}] \in {\fM}_\Gamma$, and that the
marked vertices of $\Gamma$ have distinct labels in $(X^{[3]})^T$.
The $T$-fixed points in $C_{0,1}^{(i)} \cup C_{0,2}^{(i)} \cup C_{1,2}^{(i)}$
are $Q_{i, 0},Q_{i, 1},Q_{i, 2}$.
For $0 \le j < k \le 2$, let ${\mathcal T}_{d,i, j, k}'$
be the subset of ${\mathcal T}_{d,i}'$ consisting of all
$\Gamma \in {\mathcal T}_{d,i}'$ such that the labeling $\mathfrak L$
maps the marked vertices of $\Gamma$ to $\{Q_{i,j},Q_{i,k}\}$.
Then, ${\mathcal T}_{d,i, 0, 1}'$, ${\mathcal T}_{d,i, 0, 2}'$ and
${\mathcal T}_{d,i, 1, 2}'$ form a partition of ${\mathcal T}_{d,i}'$. So
\begin{eqnarray} \label{T_di'}
\sum_{\Gamma \in \mathcal T_{d,i}'} =
\sum_{\Gamma \in \mathcal T_{d,i, 0,1}'}
+ \sum_{\Gamma \in \mathcal T_{d,i, 0, 2}'}
+ \sum_{\Gamma \in \mathcal T_{d,i, 1, 2}'}.
\end{eqnarray}

Put ${\mathcal T}_{d,i, j, k}'' = {\mathcal T}_{d,i, j, k}'/\sim$
where the relation $\sim$ is defined the same way
as in the first paragraph of subsection~\ref{subsect_Sdij'}.
%$\Gamma_1 \sim \Gamma_2$ if one can be obtained from the other
%by switching the two marked vertices.
As in (\ref{Gamma_12}) and (\ref{edij}), we get
\begin{eqnarray} \label{T_dijk'}
%   T_{d,i,j, k}' \,\,
%{\buildrel\text{def}\over=} \,\,
\sum_{\Gamma \in {\mathcal T}_{d,i, j, k}'} \int_{[M_\Gamma]^{\rm vir}}
\frac{(A \otimes B)_\Gamma - (1_X \otimes AB)_\Gamma}
{|{\bf A}_\Gamma| \,\, e(N_\Gamma^{\rm vir})}
=\gamma_{i,j,k} \cdot f_{d,i,j,k}
\end{eqnarray}
where $\gamma_{i,j,k} = -(A|_{Q_{i,j}} - A|_{Q_{i,k}})
(B|_{Q_{i,j}} - B|_{Q_{i,k}})$ and
\begin{eqnarray}  \label{f_dijk}
f_{d,i,j,k} = \sum_{\Gamma \in {\mathcal T}_{d,i,j,k}''}
\int_{[M_{\Gamma_1}]^{\rm vir}} \frac{1}{|{\bf A}_{\Gamma_1}|
\,\, e(N_{\Gamma_1}^{\rm vir})}.
\end{eqnarray}
%Likewise, let $\alpha_{i,j} = A|_{Q_{i,j}}$,
%$\beta_{i,j} = B|_{Q_{i,j}}$ and set $\gamma_{i,j,k}
%= -(\alpha_{i,j} - \alpha_{i,k})(\beta_{i,j} - \beta_{i,k})$.
By Lemma \ref{lem.fixedpoints},
we have $\gamma_{i,0,1} = -3g_{i}(w_i^2 + 2w_iz_i - 8z_i^2)^2$,
$\gamma_{i,0,2} = -3g_{i}(-8w_i^2 + 2w_iz_i  + z_i^2)^2$
and $\gamma_{i,1,2} = -243g_{i}(w_i^2 - z_i^2)^2$.
Combining (\ref{T_di'}) and (\ref{T_dijk'}) yields
\begin{eqnarray} \label{sum_T_di'}
   T_{d,i}'
&=&\sum_{\Gamma \in {\mathcal T}_{d,i}'} \int_{[\fM_\Gamma]^{\rm vir}}
   \frac{(A \otimes B)_\Gamma - (1_X \otimes AB)_\Gamma}
   {e(N_\Gamma^{\rm vir})}  \nonumber  \\
&=&\gamma_{i,0,1} \cdot f_{d,i,0,1} + \gamma_{i,0,2} \cdot f_{d,i,0,2}
   + \gamma_{i,1,2} \cdot f_{d,i,1,2}.
\end{eqnarray}

The $f_{d,i,j,k}$ can be calculated via graph sums in a manner
similar to the calculation of the $e_{d,i,j}$
in subsection \ref{subsect_Sdij'}. Note that if $f_{d,i,0,1}$ is written
as a function of the variables $w_i$ and $z_i$,
then $f_{d,i,0,2}$ can be obtained from $f_{d,i,0,1}$
by switching $w_i$ and $z_i$.
Also, for an edge $e$ of a stable graph $\Gamma$ and
for $0 \le j < k \le 2$, define $e \in [Q_{i,j}Q_{i,k}]$ if
the labeling $\mathfrak L$ of $T$ maps the two vertices of $e$ to
the set $\{ Q_{i,j}, Q_{i,k} \}$. By Lemma~\ref{15_curve},
the curves $C_{0,1}^{(i)}$, $C_{0,2}^{(i)}$
and $C_{1,2}^{(i)}$ are homologous to
$\beta_3$, $\beta_3$ and $3\beta_3$ respectively.
Therefore, for each stable graph $\Gamma$, the edges $e$ satisfy
\begin{eqnarray} \label{edge_degree}
\sum_{e \in [Q_{i, 0}Q_{i,1}]} d_e +\sum_{e \in [Q_{i, 0}Q_{i,2}]} d_e
+ \sum_{e \in [Q_{i, 1}Q_{i,2}]}3d_e = d.
\end{eqnarray}

\subsection{Cases when $1 \le d \le 4$}
\par
$\,$

When the degree $d$ is small, we can use Mathematica and
the setups of subsections~\ref{subsect_Sdij'} and \ref{subsect_Tdi'}
to make explicite computations. We now do this for $1 \le d \le 4$.

When $1 \le d \le 4$, we have verified via Mathematica that
\begin{eqnarray} \label{Sdij'}
e_{d,i,j} = \frac{w_i + z_i}{d w_iw_j(w_i - z_i)^2z_iz_j}
\quad \text{\rm and} \quad
S_{d,i,j}' = \frac{(2g_{i} + g_{j})(w_i + z_i)^3}{d w_iw_jz_iz_j}.
\end{eqnarray}
Unfortunately, we are not able to prove this formula for general $d$.

Also, for $1 \le d \le 4$, the functions $f_{d,i, 0, 1}$ are given by
\begin{eqnarray}
f_{1,i, 0,1} & = &\frac{w_i + z_i}{w_i(w_i - 2z_i)^2(w_i - z_i)z_i^2}
    \label{f_1i01}  \\
f_{2,i, 0,1} & = &\frac{2w_i^2 + 7w_i z_i + 5z_i^2}{2w_i
   (w_i - 2z_i)^2(w_i - z_i) (2w_i - z_i)z_i^2} \nonumber \\
   & = & \frac{1}{2}f_{1,i,0,1} + \frac{3(w_i + z_i)}{w_i
   (w_i - 2z_i)^2(w_i - z_i) (2w_i - z_i)z_i} \label{f_2i01}  \\
f_{3,i, 0,1} & = & \frac{2(w_i + z_i)(w_i + 4z_i)}{3w_i
   (w_i - 2z_i)^2(w_i - z_i) (2w_i - z_i)z_i^2} \nonumber \\
   & = & \frac{1}{3}f_{1,i,0,1} + \frac{3(w_i + z_i)}{w_i
   (w_i - 2z_i)^2(w_i - z_i) (2w_i - z_i)z_i} \label{f_3i01}  \\
f_{4,i,0,1} & = & \frac{2w_i^2 + 7w_i z_i + 5z_i^2}{4w_i
   (w_i - 2z_i)^2(w_i - z_i) (2w_i - z_i)z_i^2} \nonumber \\
   & = & \frac{1}{4}f_{1,i,0,1} + \frac{3(w_i + z_i)}{2w_i
   (w_i - 2z_i)^2(w_i - z_i) (2w_i - z_i)z_i}  \label{f_4i01}
\end{eqnarray}
Recall that if we regard $f_{d,i, 0,1}$ as
a function of $z_i$ and $w_i$, then $f_{d,i,0,2}$ can
be obtained from $f_{d,i, 0,1}$ by switching $z_i$ and $w_i$.
So $f_{d,i, 0,2}$ is known for $1 \le d \le 4$. Furthermore,
\begin{eqnarray}
f_{1,i,1,2} & = & 0   \label{f_1i12}  \\
f_{2,i,1,2} & = &
     \frac{w_i + z_i}{w_i(w_i - 2z_i)(w_i - z_i)^2(2w_i - z_i)z_i}
     \label{f_2i12}  \\
f_{3,i,1,2} & = &
     \frac{w_i + z_i}{w_i(w_i - 2z_i)(w_i - z_i)^2(2w_i - z_i)z_i}
     \label{f_3i12}  \\
f_{4,i,1,2} & = &
     \frac{w_i + z_i}{2w_i(w_i - 2z_i)(w_i - z_i)^2(2w_i - z_i)z_i}.
     \label{f_4i12}
\end{eqnarray}
Combining formulas (\ref{f_1i01})-(\ref{f_4i12}) with (\ref{sum_T_di'}),
we conclude that
\begin{eqnarray}
T_{1, i}'
&=&\frac{-3g_{i}(w_i^3 - 6w_i^2z_i -6w_iz_i^2 + z_i^3)}{w_i^2z_i^2}
   \label{Ti1'}  \\
T_{2, i}'
&=&\frac{-3g_{i}(w_i^3 + 12w_i^2z_i + 12w_iz_i^2 + z_i^3)}{2w_i^2z_i^2}
   = \frac{1}{2}T_{1,i}' - \frac{27g_{i}(w_i + z_i)}{w_iz_i} \label{Ti2'}\\
T_{3, i}'
&=&\frac{-3g_{i}(w_i^3 + 21w_i^2z_i + 21w_iz_i^2 + z_i^3)}{3w_i^2z_i^2}
   = \frac{1}{3}T_{1,i}' - \frac{27g_{i}(w_i + z_i)}{w_iz_i} \label{Ti3'}\\
T_{4, i}'
&=&\frac{-3g_{i}(w_i^3 + 12w_i^2z_i + 12w_iz_i^2 + z_i^3)}{4w_i^2z_i^2}
   = \frac{1}{4}T_{1,i}' - \frac{27g_{i}(w_i + z_i)}{2w_iz_i}. \label{Ti4'}
\end{eqnarray}

In view of formulas (\ref{prime}), (\ref{Sdij'})
and (\ref{Ti1'})-(\ref{Ti4'}), we obtain
\begin{eqnarray}
\langle A, B \rangle_{0,1} & = & -81
   \label{d1}  \\
\langle A, B \rangle_{0,2} & = & -\frac{81}{2} + 81= \frac{81}{2}
   \label{d2}  \\
\langle A, B \rangle_{0,3} & = & -\frac{81}{3} + 81= 54
   \label{d3}  \\
\langle A, B \rangle_{0,4} & = & -\frac{81}{4} + \frac{81}{2}= \frac{81}{4}.
   \label{d4}
\end{eqnarray}
%A natural guess for the general form for $d > 1$ is
%\begin{equation}
%\langle A, B \rangle_d  = -\frac{81}{d} + 81.
%\end{equation}

\begin{proposition} \label{4_prop}
Let $X = \PP^2$, and $\ell \subset X$ be a line.
Then, the $2$-point genus-$0$ Gromov-Witten invariant
$\langle \pd(\frak a_{-3}(\ell)\vac),
\pd(\frak a_{-3}(X)\vac) \rangle_{0, d}$ is equal to
$-27$, $27/2$, $18$ and $27/4$
when $d$ is equal to $1$, $2$, $3$ and $4$ respectively.
\end{proposition}
\begin{proof}
Follows immediately from (\ref{tilde_A_to_A}) and (\ref{d1})-(\ref{d4}).
\end{proof}

%\begin{remark}
%The Gromov-Witten invariant $\langle A, B \rangle_{0,d}$ appears
%to be given by a formula of the form $\frac{-81}{d} + \frac{81}{f(d)}$ for
%$d\geq 2$, where $f(d)$ is some unknown function. Since $f(2) = f(3) = 1,
%f(4) =2$ one possible guess is that $f(d)$ is the $(d-1)$'st Fibonacci number.
%However, we have no evidence for this guess other than the first 3 terms.
%\end{remark}

%%
%%
%%
%%
%%
%%
%%
%%
%%

\end{document}